\documentclass[A4paper,reqno]{amsart}
\usepackage[usenames,dvipsnames]{xcolor}
\usepackage[T1]{fontenc}
\usepackage[utf8]{inputenc}
\usepackage{lmodern}
\usepackage[UKenglish]{babel}
\usepackage{amsmath, amssymb, amsfonts, amsthm, mathtools, bbm}
\usepackage{extarrows}
\usepackage{mathrsfs}
\usepackage{pgfplots}
\usetikzlibrary{arrows}
\usepackage{tikz-cd}
\usepackage{graphicx}
\usepackage{faktor}
\usepackage{young}
\usepackage[centertableaux]{ytableau}
\usepackage{verbatim} 
\usepackage{comment}
\usepackage{wasysym}
\usepackage{enumitem}
\usepackage{microtype}
\usepackage{csquotes}
\usepackage{cancel}
\usepackage{color}      
\usepackage{xparse}
\usepackage{longtable}
\usepackage{hyperref}
\usepackage{cleveref}	
\usepackage{braids}
\usepackage{adjustbox}
\usepackage{forest} 
\usepackage{dynkin-diagrams}
\usepackage{rank-2-roots}
\definecolor{light-gray}{gray}{0.85}

\DeclareMathOperator{\coireal}{\mathcal{M}} 
\DeclareMathOperator{\coidiag}{\mathscr{M}} 
\DeclareMathOperator{\algreal}{\mathcal{V}} 
\DeclareMathOperator{\algdiag}{\mathscr{V}} 
\binoppenalty=\maxdimen         
\relpenalty=\maxdimen           

\newcommand{\opn}[1]{\operatorname{#1}}

\newtheorem{theo}{Theorem}
\numberwithin{theo}{section}
\newtheorem{lemm}[theo]{Lemma}
\newtheorem{prop}[theo]{Proposition}
\newtheorem{koro}[theo]{Corollary}

\theoremstyle{definition}
\newtheorem{defi}[theo]{Definition}
\newtheorem{beis}[theo]{Example}
\newtheorem{rema}[theo]{Remark}

\DeclareMathOperator{\Hom}{Hom}
\DeclareMathOperator{\End}{End}

\DeclareMathOperator{\C}{\mathbb{C}}

\DeclareMathOperator{\sgn}{sgn}
\DeclareMathOperator{\id}{id}

\newcommand{\bp}{b_{n,{+}, \varepsilon}}
\newcommand{\bn}{b_{n,{-}, \varepsilon}}

\DeclareMathOperator{\modd}{\mathbf{Rep}}

\DeclareMathOperator{\TL}{TL}

\newcommand{\dotTL}{\opn{TL}(B_n)}
\DeclareMathOperator{\Ug}{U}
\DeclareMathOperator{\sl2}{\mathfrak{sl}_2}

\DeclareMathOperator{\lie}{\mathfrak{g}}
\DeclareMathOperator{\gl}{\mathfrak{gl}}

\newcommand{\Uq}{\Ug_q(\sl2)}

\newcommand{\Vq}{\Ug_q'(\gl_1)}

\DeclareMathOperator{\B}{B}

\DeclareMathOperator{\mZ}{\mathbb{Z}}

\DeclareMathOperator{\mN}{\mathbb{N}_0}
\DeclareMathOperator{\mU}{\mathbbm{1}}

\newcommand{\Hecke}{\mathcal{H}}

\newcommand{\shift}[1]{\langle#1\rangle} 

\DeclareMathOperator*{\tlcupcap}{\,\pretlcupcap\,}

\DeclareMathOperator*{\tlline}{\pretlline}
 
\newcommand{\cD}{{\mathcal{D}}} 

\usetikzlibrary{calc}
\usetikzlibrary{tqft}
\newcommand{\cbox}[1]{ \vcenter{ \hbox{#1} } }
\usetikzlibrary{patterns,decorations.pathreplacing}
\tikzset{
	tldiagram/.style={thick, scale=0.35}
}
\newcommand{\tlcoord}[2]
{
	(2*#2 , 3*#1)
}

\newcommand{\lineup}{-- ++(0,3)}
\newcommand{\halflineup}{-- ++(0,1.5)}
\newcommand{\quarterlineup}{-- ++(0,0.75)}
\newcommand{\linedown}{-- ++(0,-3)}
\newcommand{\halflinedown}{-- ++(0,-1.5)}

\newcommand{\smalllineup}{-- ++(0,2)}
\newcommand{\smalllinedown}{-- ++(0,-2)}
\newcommand{\linewave}[2]{
	.. controls +(0,1.5*#1) and +(0,1.5*-#1) .. ++(2*#2, 3*#1)
}

\newcommand{\dlineup}{\halflineup \onedot \halflineup}
\newcommand{\halfdlineup}{\quarterlineup \onedot \quarterlineup}

\newcommand{\capright}{arc (180:0:1)}
\newcommand{\capleft}{arc (0:180:1)}
\newcommand{\cupright}{arc (180:360:1)}
\newcommand{\cupleft}{arc (360:180:1)}

\newcommand{\onedot}{node {$\bullet$}}

\newcommand{\dcapright}{arc (180:90:1) node {$\bullet$} arc (90:0:1)}

\newcommand{\dcupright}{arc (180:270:1) node {$\bullet$} arc (270:360:1)}

\newcommand{\dcaprightsmall}{arc (180:90:1) node {\scalebox{0.7}{\textbullet}} arc (90:0:1)}

\newcommand{\dxcapright}[1]{arc (180:90:#1) node {$\bullet$} arc (90:0:#1)}

\newcommand{\dcuprightsmall}{arc (180:270:1) node {\scalebox{0.7}{\textbullet}} arc (270:360:1)}

\newcommand{\xcapright}[1]{arc (180:0:#1)}

\newcommand{\xcupright}[1]{arc (180:360:#1)}
\newcommand{\xcupleft}[1]{arc (360:180:#1)}

\newcommand{\maketlboxnormal}[2]{ 
	++(-1,-0.75)
	+(-0.5, -0.3)
	[fill=white!20!white, rounded corners]
	rectangle
	+(#1*2 + 0.5, 1.5 + 0.3)
	node at ++(#1,0.75) {#2}
}
\newcommand{\maketlboxred}[2]{  
	++(-1,-0.75)
	+(-0.5, -0.3)
	[fill=red!20!white, rounded corners]
	rectangle
	+(#1*2 + 0.5, 1.5 + 0.3)
	node at ++(#1,0.75) {#2}
}
\newcommand{\maketlboxgreen}[2]{
	++(-1,-0.75)
	+(-0.5, -0.3)
	[fill=ForestGreen!20!white, rounded corners]
	rectangle
	+(#1*2 + 0.5, 1.5 + 0.3)
	node at ++(#1,0.75) {#2}
}
\newcommand{\maketlboxblue}[2]{ 
	++(-1,-0.75)
	+(-0.5, -0.3)
	[fill=cyan!20!white, rounded corners]
	rectangle
	+(#1*2 + 0.5, 1.5 + 0.3)
	node at ++(#1,0.75) {#2}
}

\newcommand{\negtwista}[2]{
	\draw [ultra thick] \tlcoord{0+3*#1}{0+2*#2} -- \tlcoord{0.45+3*#1}{0.4+2*#2};
	\draw [ultra thick] \tlcoord{0.55+3*#1}{0.6+2*#2} -- \tlcoord{1+3*#1}{1+2*#2};
	\draw  \tlcoord{0+3*#1}{1+2*#2} -- \tlcoord{1+3*#1}{0+2*#2};
}

\newcommand{\negtwist}[2]{
	\negtwista{#1}{#2}
	\negtwist{1+#1}{#2}
}

\newcommand{\pretlcupcap}
{
	\begin{tikzpicture}[scale=0.092]
		\draw \tlcoord{1}{0} \cupright;
		\draw \tlcoord{0}{0} \capright;
	\end{tikzpicture}
}
\newcommand{\predtlcupcap} 
{
	\begin{tikzpicture}[scale=0.092]
		\draw \tlcoord{1}{0} \dcuprightsmall;
		\draw \tlcoord{0}{0} \dcaprightsmall;
	\end{tikzpicture}
}
\newcommand{\tlcircle}
{
	\begin{tikzpicture}[scale=0.11]
		\draw (0,0) circle (1);
	\end{tikzpicture}
}
\newcommand{\pretlline}
{
	\vert
}

\newcommand{\blueminus}{\textcolor{blue}{-}}
\newcommand{\redbullet}{\textcolor{red}{\bullet}}
\newcommand{\posredbluecrossing}[2]{
	\draw[red] \tlcoord{#1}{#2} \linewave{1}{1};
	\draw[white, double=blue] \tlcoord{#1}{#2+2} \linewave{1}{-1} ;
}
\newcommand{\posredredcrossing}[2]{
	
	\draw[red] \tlcoord{#1}{#2+2} \linewave{1}{-1} ;
	\draw[white, double=red] \tlcoord{#1}{#2} \linewave{1}{1};
}
\newcommand{\posbluebluecrossing}[2]{
	
	\draw[blue] \tlcoord{#1}{#2+2} \linewave{1}{-1} ;
	\draw[white, double=blue] \tlcoord{#1}{#2} \linewave{1}{1};
}
\newcommand{\posblueredcrossing}[2]{
	\draw[red] \tlcoord{#1}{#2+2} \linewave{1}{-1} ;
	\draw[white, double=blue] \tlcoord{#1}{#2} \linewave{1}{1};
}

\newcommand{\littlebluered}{
	\!\!\cbox{
		\begin{tikzpicture}[tldiagram,scale=3/10]
			\posblueredcrossing{0}{0}
		\end{tikzpicture}
	}\!\!
}
\newcommand{\littleredblue}{\!\!\cbox{
		\begin{tikzpicture}[tldiagram,scale=3/10]
			\posredbluecrossing{0}{0}
		\end{tikzpicture}
	}\!\!}

\newcommand{\widerhakenup}{ -- ++ (0.4,0) -- ++(0,-0.5) -- ++(0,0.5) -- ++ (-0.8,0) -- ++(0,-0.5) -- ++(0,0.5) -- ++ (0.4,0)
}
\newcommand{\widerhakendown}{-- ++ (0.4,0) -- ++(0,0.5) -- ++(0,-0.5) -- ++ (-0.8,0) -- ++(0,0.5) -- ++(0,-0.5) -- ++ (0.4,0)}
\newcommand{\bluewiderhakenup}{\!\!
	\cbox{
		\begin{tikzpicture}[tldiagram, scale=5/10]
			\draw[blue] \tlcoord{0}{0} \halflineup \widerhakenup; 
		\end{tikzpicture}
	}\!\!
}
\newcommand{\bluewiderhakendown}{\!\!
	\cbox{
		\begin{tikzpicture}[tldiagram, scale=5/10]
			\draw[blue] \tlcoord{0}{0} \halflinedown \widerhakendown; 
		\end{tikzpicture}
	}\!\!
}

\setlist[enumerate]{leftmargin=*, label= \roman*)}

\setlist{
	listparindent=\parindent,
	parsep=0pt
}

\begin{document}
\title[Smallest quantum coideal and Jones--Wenzl projectors]{{Diagrammatics for the smallest quantum coideal and Jones--Wenzl projectors}}
\author[C. Stroppel]{Catharina Stroppel}
\address{C. S.: Department of Mathematics, University of Bonn, Germany}
\email{stroppel@math.uni-bonn.de}
\author[Z. Wojciechowski]{Zbigniew Wojciechowski}
\address{Z. W.: Department of Mathematics, University of Dresden, Germany}
\email{zbigniew.wojciechowski@tu-dresden.de}

		\maketitle
		\begin{abstract}
			We describe algebraically, diagrammatically and in terms of weight vectors, the restriction of tensor powers of the standard representation of quantum $\sl2$  to a coideal subalgebra.  We realise  the category as module category over the monoidal category of type $\pm1$ representations in terms of string diagrams and via generators and relations. The idempotents projecting onto the quantized eigenspaces are described as type $B/D$ analogues of Jones--Wenzl projectors. As an application we introduce and give recursive formulas for analogues of $\Theta$-networks.
		\end{abstract}
		\section{Introduction}
				
		Decomposing tensor products of representations for a group $G$, a Lie (super-) algebra $\mathfrak{g}$ and similar structures into indecomposable summands is key for understanding representations of the given structure. 
		In the world of complex semisimple Lie algebras, the prime example is $\mathfrak{sl}_2$, whose tensor products are described in terms of the Clebsch--Gordan rule. This description lifts to the quantum setting when replacing $\sl2$ with $\Uq$. 
		In the world of quantum symmetric pairs (QSP's) the prime example is $\Ug_q'(\sl2^{ \!\iota})\subset \Uq$. Loosely speaking a QSP is a quantization of the canonical inclusion $\Ug(\lie^{ \iota})\subset \Ug(\lie)$ of universal enveloping algebras coming from an embedding of Liealgebras $\lie^{\iota}\subset \lie$ given by a classical symmetric pair (that is an embedding of the Lie algebra of fixed points $\lie^{\iota}$ under a Lie algebra involution $\iota\colon \lie\to \lie$). The resulting subalgebra $\Ug_q'(\lie^{\iota})\subset \Ug_q(\lie)$ is not a Hopf subalgebra in general, but just a one-sided coideal. This has consequences for the representations of $\Ug_q'(\lie^{\iota})$ and their morphisms: they cannot be tensored together; instead, representations of $\Ug_q(\lie)$ \emph{act} on this category. More precisely, the category of representations of  $\Ug_q'(\lie^{\iota})$ is a (even braided) module category over the (braided) monoidal category of representations of $\Ug_q(\lie)$, in the sense of \cite[\S5.1]{brochier23}, \cite[\S3 Def.7]{HaOl-actions-tensor-categories}. The purpose of this paper is to make this totally explicit in case of  $\Ug_q'(\sl2^{ \!\iota})\subset \Uq$ using diagrammatic methods.  
		
We extend a well-established connection between representation theory and diagrammatic from the familiar type $A$ to other Lie types by developing a generalisation of the Temperley-Lieb recoupling theory,\cite{flath95}, \cite{kauffman1994}. The paper partially presents results obtained in the thesis  \cite{Woj2022}, to which we also refer for more background material. Since the original theory in type $A$ is ultimately based on quantum Schur--Weyl duality or quantum Skew-Howe duality \cite{cautis14}, the starting point of this paper are the Schur--Weyl and skew Howe duality results between a specialized 2-parameter type $B$ Hecke algebra (which contains the type $D$ Hecke algebra) and the type AIII quantum symmetric pair, \cite{bao18}, \cite{stroppel2018}. These results and the invention and categorification of $i$-canonical bases in \cite{bao18}, \cite{stroppel2018}, suggest that the representation theory of the smallest coideal $\Vq\subset \Uq$ can be studied using Temperley--Lieb categories of type $D$ and described by the corresponding blob-diagrammatics from \cite{green1998}, \cite{tomdieck1994}, \cite{tomdieck1998}. 
Such generalized Temperley--Lieb categories appeared in connection with $\Uq$ already in \cite{iohara18} to describe intertwiners between a dominant Verma tensored with copies of the standard representation. Despite the fact that the coideal $\Vq$ was already a prominent example in the theory of quantum homogeneous spaces, \cite{mueller99}, before quantum symmetric pairs were classified in  \cite{letzter02}, its representation category has not been studied in detail yet, but see remarkable first steps see e.g. \cite[Prop~3.14]{watanabe2021} (a study of weight modules and classification results) and  \cite{Watanabecrys} (for crystal bases).

\subsection*{Overview and results} We start our study of tensor products and fusion rules in \Cref{section on fusion rules} in terms of deformed weight vectors. These are interesting eigenvectors for a commuting subalgebra which (in contrast to the quantum group case) does not consists of primitive elements. This means in particular that tensor products of weight vectors are not necessarily weight vectors.  
Various normalisations problems of the Jones polynomial and other invariants related to $\Uq$ encode subtle sign issues, see e.g. \cite{tingley2015minus}, \cite{StICM} which can be explained by the interplay of type $1$ and type $-1$ representations. In \Cref{section string calculus} we therefore introduce simultaneous diagrammatics for type $1$ and $-1$ representations of $\Uq$ and describe diagrammatically the module category of representations of $\Vq$, which arise as restrictions of tensor products of representations of $\Uq$. In \Cref{section on type B projectors} we relate this to the decompositions from \Cref{section on fusion rules}. We define analogues of Jones--Wenzl projectors, prove uniqueness properties and give a new combinatorial method to count dimensions of morphisms spaces for the coideal. We end the paper in \Cref{section application to theta networks} by  introducing and then calculating certain $\Theta$-networks of type $D$. 

All results of this paper which are formulated in terms of Temperley--Lieb categories, in particular the Temperley--Lieb algebra of type $D$, can be categorified easily (and straight-forward) in terms of (geometric bi)modules over the type $D$ Khovanov arc algebras  \cite{stroppel2016}, that is a type $D$ version of Khovanov's famous arc algebra which he used to categorify the Reshethikhin--Turaev $\Uq$-invariant of tangles. We however do not pursue this here. Instead we explicitly write out the decategorified picture in a way such that every step can be categorified easily. In type $A$, the Khovanov arc algebra can be realised, see see \cite{BS3},  as a cyclotomic quotient of a KLR-algebra \cite{Mathas}. These algebras, introduced in \cite{KL}, \cite{R}, are generalisations of Nil-Hecke algebras which can be used to category the positive part of $\Uq$. In parallel, the coideal $\Vq$ has been categorified recently in \cite{brundan2023} via an KLR approach using projective modules over the Nil--Brauer category. A direct connection between the two possible categorifications is not known at the moment and will be subject of a further study. 
		\subsubsection*{Acknowledgments}
		The second author thanks Gaetan Borot and the Humboldt-Universität zu Berlin for the CIMPA school in Puerto Madryn in 2023. Both authors were supported by the Hausdorff Center of Mathematics (EXC 2047) in Bonn. 
			\section{Fusion rules and decompositions for the coideal} \label{section on fusion rules}
		Fix as ground field $k=\C(q)$. For $n\in \mZ$ set $[n]=\frac{q^n-q^{-n}}{q-q^{-1}}\in k$. 
		\begin{defi}
			The \textit{smallest quantum group} $\Uq$ is the $\C(q)$-algebra generated by $E,F,K,K^{-1}$ subject to the relations
			\[
			KE=q^2EK, \, KF=q^{-2}KF, \, EF-FE=\frac{K-K^{-1}}{q-q^{-1}},  \, KK^{-1}=1=K^{-1}K.
			\]
			It is a Hopf algebra with comultiplication given by
			\[
			\Delta(E)=E\otimes 1 + K \otimes E, \quad \Delta(F)=F\otimes K^{-1} + 1 \otimes F, \quad \Delta(K)=K\otimes K.
			\]
		\end{defi}
	\begin{defi} \label{definition standard representation} \label{definition trivial representation}
		Let $m\geq 0$, $\varepsilon\in\{{1},{-1}\}$. The \textit{$m+1$-dimensional irreducible representation $V_{m,\varepsilon}$ of type $\varepsilon$} is the $m+1$-dimensional $k$-vector space with basis $x_{0,\varepsilon},\ldots, x_{m,\varepsilon}$ and $\Uq$-action
		\begin{gather*}
			Ex_{i,\varepsilon}=\begin{cases}
				[m-i+1]x_{i-1,\varepsilon} & i>0, \\
				0 	& i=0,
			\end{cases} \quad Fx_{i,\varepsilon}=\begin{cases}
				\varepsilon[i+1]Fx_{i+1,\varepsilon} & i>0, \\
				0 	& i=m,
			\end{cases} \\
			Kx_{i,\varepsilon}=\varepsilon q^{m-2i}x_{i,\varepsilon}.
		\end{gather*}
		Set $V_{\varepsilon}\coloneqq V_{1,\varepsilon}$ the \textit{type $\varepsilon$ standard representation}, and  $k_{\varepsilon}\coloneqq V_{0,\varepsilon}$ the \textit{type $\varepsilon$ trivial} representation. We write $x_{\varepsilon}\coloneqq x_{0,\varepsilon}, y_{\varepsilon}\coloneqq x_{1,\varepsilon}$ and $1_{\varepsilon}\coloneqq x_{0,\varepsilon}$ for their bases.
	\end{defi}
		Consider $B\coloneqq q^{-1}EK^{-1}+F\in \Uq$. We directly obtain the following:
		\begin{lemm} \label{lemma that we have coideal}
			The element $B\in \Vq$ satisfies $\Delta(B)=B\otimes K^{-1} + 1 \otimes B$.
		\end{lemm}
		\begin{defi}
			The \textit{smallest type AIII right coideal} $\Vq$ is the subalgebra of $\Uq$ generated by the element $B\coloneqq q^{-1}EK^{-1}+F$.
		\end{defi}
		\begin{rema}
			By \Cref{lemma that we have coideal}, $\Vq$ is a right coideal subalgebra of $\Uq$, justifying the name. The coideal subalgebra $\Vq$ of $\Uq$ appears as one of the fundamental examples of quantum symmetric pairs \cite[\S7]{letzter03}.
		The left coideal version $'\!\Ug_q(\gl_1)$ of $\Vq$ is generated by $'\!\B\coloneqq E+qKF$. 
		\end{rema}
		The following describes tensor powers of $V_{\varepsilon}$ over $\Vq$.
		\begin{theo}[Coideal tensor product decomposition] \label{Theorem decomposition of tensor power of V}
			Let $\varepsilon\in\{1,-1\}$. Consider for $n\in \mZ$ the special elements in $V_{\varepsilon}$ given by
			\[
			v_{n, \varepsilon}= x_{\varepsilon} + \varepsilon q^{\varepsilon n}y_{\varepsilon}, \quad w_{n, \varepsilon}=x_{\varepsilon}-\varepsilon q^{-\varepsilon n}y_{\varepsilon}.
			\]
			\begin{enumerate}
				\item \label{part i of gl1 x gl1 decomp} Then $V_{\varepsilon}$ decomposes over $\Vq$ into summands
				\[
				V = L(1)\oplus L(-1)\coloneqq \shift{v_{0, \varepsilon}} \oplus \shift{w_{0, \varepsilon}},
				\]
				which are eigenspaces for $B$ with eigenvalues $1$ respectively $-1$.
				\item \label{part ii of gl1 x gl1 decomp} Let $z$ be an eigenvector for $B$ with eigenvalue $[n]$, $n\in\mZ$. Then $z\otimes v_{n, \varepsilon}$ and $z\otimes w_{n,\varepsilon}$ are eigenvectors for $B$ with eigenvalue $[\varepsilon n+1]$ and $[\varepsilon n-1]$ respectively.
				\item \label{part iv of gl1 x gl1 decomp} The $n$-th tensor power of $V_{\varepsilon}$ decomposes into $1$-dimensional subrepresentations $L([n-2k])$ which are eigenvectors for $B$ with eigenvalue $[n-2k]$,
				\begin{equation} \label{tensor decomposition}
					V_{\varepsilon}^{\otimes n} = \bigoplus_{k=0}^{n} {n \choose k} L([n-2k]).
				\end{equation}
			\end{enumerate}
		\end{theo}
		\begin{proof} 
				Part \ref{part i of gl1 x gl1 decomp} is clear. With \Cref{lemma that we have coideal} we calculate
				\begin{align*}
					B (z\otimes (x_{\varepsilon}\pm \varepsilon q^{\pm \varepsilon n}y_{\varepsilon}))
					&=Bz \otimes K^{-1}(x_{\varepsilon}\pm \varepsilon q^{\pm \varepsilon n}y_{\varepsilon}) + z \otimes B(x_{\varepsilon}\pm\varepsilon q^{\pm \varepsilon n} y_{\varepsilon}) \\
					&= [n] z \otimes (\varepsilon q^{-1} x_{\varepsilon} \pm q^{\pm\varepsilon n+1}y_{\varepsilon}) + z \otimes (\varepsilon y_{\varepsilon} \pm q^{\pm \varepsilon n} x_{\varepsilon}) \\
					&=[\varepsilon n\pm1] z \otimes (x_{\varepsilon}\pm \varepsilon q^{\pm \varepsilon n}y_{\varepsilon}),
				\end{align*}
				\[
				\text{since }
				\varepsilon q^{-1}[n] \pm  q^{\pm \varepsilon n} =[\varepsilon n\pm 1], \quad  q^{\varepsilon n+1}[\pm n]+\varepsilon=\pm \varepsilon q^{\pm \varepsilon n} [\varepsilon n\pm1].
				\]
		This shows \ref{part ii of gl1 x gl1 decomp}. Part \ref{part iv of gl1 x gl1 decomp} follows from \ref{part i of gl1 x gl1 decomp} and \ref{part ii of gl1 x gl1 decomp}.
		\end{proof}
		\begin{koro}[Clebsch--Gordan rule] \label{general characters}
			For $m\geq 0$ denote by $V_{m,\varepsilon}$ the $m+1$-dimensional irreducible representation of $\Uq$ of type $\varepsilon$. Then
			\[
			L([n])\otimes V_{m,\varepsilon}=\bigoplus_{k=0}^mL([n+m-2k]).
			\]
		\end{koro}
		\begin{proof}
			This follows from \Cref{Theorem decomposition of tensor power of V} using $V_{n,\varepsilon}\otimes V_{1}\cong V_{n+1,\varepsilon}\otimes V_{n-1,\varepsilon}$. 
		\end{proof}

	\begin{rema}
		For $\varepsilon=1$ and $n=0$,  \Cref{general characters} recovers the decomposition of \cite[\S4]{mueller99}. Concretely consider the element $x=qEK^{-1}+F$, obtained by choosing $\alpha=0, \beta=q, \gamma=1$ in their setup. Then \cite[Proposition 4.2]{mueller99} states that the eigenvalues of the right action of $x$ on $V_{m,1}$ are of the form $\{[m-2k]\vert k\in\{0,\ldots, m\}$. Our element $B$ is given by $B=f(S'(-x))$, where $S'$ is the antipode with respect to their comultiplication (call it $\Delta'$) and $f$ is the $\C$-linear Hopf algebra isomorphism
		\begin{eqnarray*}
			&f\colon \Uq\to \Uq, \quad
			f(E)= EK^{-1}, f(F)=KF, f(K)=K^{-1}, F(q)= q^{-1}.&
		\end{eqnarray*}
	\end{rema}
		\begin{beis} \label{beis decomposition}
			By \Cref{general characters}, $(L([2])\oplus L(-[2])) \otimes V_{2,1}\otimes V_{2,1}$ decomposes as 
			\[
			L([6])\oplus L([4])^2\oplus L([2])^4 \oplus L(0)^4 \oplus L([-2])^4 \oplus L([-4])^2 \oplus L([-6]).
			\]
			In particular $L(0)$ appears with multiplicity $4$, which implies by adjunction
			\[
				\dim \Hom_{\Vq}((L([2])\oplus L([-2]))\otimes V_{2,1}, V_{2,1})=4.
			\] 
		\end{beis}
		\begin{rema}
			The coideal $\Vq\subset\Uq$ can be extended to a coideal $\Ug_q'(\gl_1\times \gl_1)\subset \Ug_q(\gl_2)$. Concretely if $D_1, D_{-1}$ are the additional group-like generators of $\Ug_q(\gl_2)$ subject to the additional relations
			\[
				D_iE=q^{i}E_iD_i, \, D_iF=q^{-i}FD_i, \,  D_1D_{-1}=D_{-1}D_1, \, K=D_1D_{-1}^{-1}
			\]
			we let $\Ug_q'(\gl_1\times \gl_1)$ be the subalgebra generated by $B$ and $D\coloneqq D_1D_2$.  This coideal has the advantage that $D$ distinguishes different tensor powers of $V_{\varepsilon}$. 
		\end{rema}
		\section{Different types of Temperley--Lieb algebras} \label{section temperley-lieb}
		We describe the endomorphism algebra of $V_{\varepsilon}^{\otimes n}$ using a specialized type $B$ Temperley--Lieb algebra $\TL(B_n)$. We fix $\varepsilon\in \{1,-1\}$.
		\begin{defi} \label{Definition Temperley--Lieb of Type B}
			Let $n\geq1$. The (specialized) \textit{type $B$ Temperley--Lieb algebra $\TL(B_n)_{\varepsilon}$} is the $\C(q)$-algebra with generators $s_0, U_1, \ldots, U_{n-1}$ which are subject to the \emph{type $A$ Temperley--Lieb relations} for $i,j\in\{1,\ldots, n-1\}$, $|i-j| \geq 2$:
			\begin{equation}
			U_i^2=-\varepsilon(q+q^{-1})U_i,\quad U_iU_j=U_jU_i,\\
			U_i U_{i+1} U_i=U_i, U_{i+1} U_{i} U_{i+1}=U_{i+1}
			\end{equation}
			and the additional, specialized type \emph{$B$ Temperley--Lieb relations} for $i\in\{2,\ldots,n-1\}$:
				\begin{equation}
				s_0^2=1,\quad U_1 s_0  U_1=0,\quad, s_0  U_i = U_i s_0
				\end{equation}
			involving the generator $s_0$. The subalgebras generated by $U_i$ for $1\leq i\leq n-1$, respectively $U_i$ for $0\leq i \leq n$ with $U_0\coloneqq s_0U_1s_0$ are the \textit{Temperley--Lieb algebras of type $A_{n-1}$ and type $D_n$}, short $\TL(A_{n-1})_\varepsilon$ respectively $\TL(D_{n})_\varepsilon$. 
		\end{defi}
			
		\begin{rema} \label{our definition vs Greens}
			The algebra $\TL(B_n)_{\varepsilon}$ appears as $T'D_n$ in \cite[Definition 2.3]{tomdieck1998}, our $s_0$ is $\varepsilon_0-1$ in their notation. It is closely related to the algebra $TB_n$ introduced in \cite[\S4]{tomdieck1994}, where a basis theorem was established. The blob diagrammatics description for $\TL(B_n)$ and $\TL(D_n)$ was established in \cite[Thm 4.1, 4.2]{green1998}.
			Our algebra $\TL(D_n)_{\varepsilon}$ is the quotient of Green's type $D_n$ Temperley--Lieb algebra by the relation $U_0U_1=0=U_1U_0$, which means the closed dotted loop is zero. In \cite[Theorem 2.7]{tomdieck1998} the relation between all Temperley--Lieb algebras can be found, our algebra $\TL(D_n)$ is $T''(D_n)$ in the notation from there.
		\end{rema}
		The following lemma allows to define an involution $\sigma$ on $\TL(B_n)$ and gives a realization of $\TL(D_n)$ as its fixed-point subalgebra. 
		\begin{lemm} \label{Type B automorphism s_o flips sign} 
			There is a well-defined $k$-linear algebra involution
			\begin{align*}
				\sigma\colon \dotTL &\to \dotTL , \quad 
				s_0 &\mapsto -s_0\quad 
				U_i &\mapsto U_i \quad  \text{ where } i=1,\ldots,n-1.
			\end{align*} 
			The fixed-point subalgebra $\TL(B_n)^{\sigma}$ for this involution is precisely $\TL(D_n)$.
		\end{lemm}
		\begin{proof}
			One checks the assignments yield an algebra involution. The second part follows from the diagrammatic description of $\TL(D_n)$ in \cite[Thm 4.2]{green1998} as span of all type $B$ diagrams decorated with an even number of dots and \Cref{our definition vs Greens}.
		\end{proof}
		
		Next we recall the connection between $\TL(B_n)$ and a type $B$ Hecke algebra.
		
		\begin{defi} \label{definition Hecke algebra}
			The (specialized) Hecke algebra $\Hecke_{1,q}(B_n)$ of type $B_n$ is the $\C(q)$-algebra with generators $s_0, H_1, \ldots, H_{n-1}$, subject to the type $A$ Hecke relations
			\begin{enumerate}[label = {(H\arabic*)}, align=left]
				\item \label{Hecke2 parameter relation}
				$H_i^2=1 + (q^{-1}-q)H_i$ for all $1\leq i\leq n-1$,
				\item \label{Hecke usual Braid relation}
				$H_iH_{i+1}H_i=H_{i+1}H_iH_{i+1}$ for all $1\leq i \leq n-2$,
				\item \label{Hecke usual comm relation}
				$H_iH_j=H_j H_i$ for all $0\leq i,j \leq {n-1}$ such that $|i-j|\geq 2$,
			\end{enumerate}
			and the additional, specialized type $B$ Hecke relations
			\begin{enumerate}[label = {(BH\arabic*)}, align=left]
				\item \label{Hecke1 parameter relation}
				$s_0^2=1$,
				\item \label{Hecke type B relation}
				$s_0H_1s_0H_1=H_1s_0H_1s_0$,
				\item \label{Hecke comm relation for extra generator} $s_0H_i=H_is_0$ for all $2\leq i \leq n-1$.
			\end{enumerate}
		\end{defi}
		\begin{lemm}
			There is a surjective $\C(q)$-algebra homomorphism $\varphi_{\varepsilon}\colon\Hecke_{1,q}(B_n)\rightarrow \TL(B_n)_{\varepsilon}$ defined on the generators via the assignments $\varphi_{\varepsilon}(s_0)= s_0, \varphi_{\varepsilon}(H_i)= U_i + \varepsilon q^{-\varepsilon}$, where $i=1,\ldots,n-1$.
		\end{lemm}
		\begin{proof}
			This is a reformulation of \cite[Proposition 2.12]{tomdieck1998}. In the notation from \cite[Proposition 2.12]{tomdieck1998} $\Hecke_{1,q}(B_n)$ is $H'(D_n)$, and the map from there is our map precomposed with the $q$-linear automorphism sending  $s_0$ to $s_0$ and $H_i$ to $-H_i^{-1}$.
		\end{proof}
We next show that the action of $\Hecke_{1,q}(B_n)$ on $V_{\varepsilon}^{\otimes n}$ factors through $\TL(B_n)$. 
		\begin{defi}\label{endos}
			We define $\C(q)$-linear endomorphisms
			\begin{equation*}
				\kappa_{V_{\varepsilon}}\colon V_{\varepsilon}\to V_{\varepsilon}, \quad 
				x_{\varepsilon}\mapsto y_{\varepsilon}, \quad y_{\varepsilon}\mapsto x_{\varepsilon}, \quad\text{and}\quad R_{V_{\varepsilon},V_{\varepsilon}} \colon V_{\varepsilon} \otimes V_{\varepsilon} \rightarrow V_{\varepsilon} \otimes V_{\varepsilon} \text{ given by}
				\end{equation*}
				\begin{eqnarray*}
&\begin{array}{lcllcl}
x_{\varepsilon}\otimes x_{\varepsilon} &\mapsto& \varepsilon q^{-\varepsilon} x_{\varepsilon} \otimes x_{\varepsilon} &
				x_{\varepsilon} \otimes y_{\varepsilon} &\mapsto& y_{\varepsilon} \otimes x_{\varepsilon} + \delta_{\varepsilon,-1}(q^{-1}-q)x_{\varepsilon}\otimes y_{\varepsilon}\\
				y_{\varepsilon}\otimes y_{\varepsilon} &\mapsto& \varepsilon q^{-\varepsilon} y_{\varepsilon}\otimes y_{\varepsilon},&
				y_{\varepsilon} \otimes x_{\varepsilon} &\mapsto& x_{\varepsilon} \otimes y_{\varepsilon} + \delta_{\varepsilon,1}(q^{-1}-q)y_{\varepsilon}\otimes x_{\varepsilon}.
				\end{array}&
			\end{eqnarray*}
		\end{defi}
		\begin{lemm} \label{map from hecke to endos}
			There is a well-defined algebra homomorphism such that
			\begin{align*}
				\Phi_{\varepsilon}\colon\Hecke_{1,q}(B_n)&\rightarrow \End_{\C(q)}(V_{\varepsilon}^{\otimes n}) \\
				s_0 &\mapsto \kappa_{V_{\varepsilon}} \otimes \id_{V^{\otimes (n-1)}}, \\
				H_i &\mapsto \id_{V^{\otimes (i-1)}} \otimes R_{V_{\varepsilon},V_{\varepsilon}} \otimes \id_{V^{\otimes (n-i-1)}} \text{ where } i=1,\ldots,n-1.
			\end{align*}
		\end{lemm}
		\begin{proof}This (and well-definedness in \Cref{endos}) follows from the definitions.
		\end{proof}
		\begin{defi}
			We define an endomorphism $\underline{H}_{V_\varepsilon,V_\varepsilon}$ of $V_{\varepsilon}\otimes V_{\varepsilon}$ by 
			\[
			\underline{H}_{V_\varepsilon,V_\varepsilon}\coloneqq R_{V_{\varepsilon},V_{\varepsilon}} - \varepsilon q^{-\varepsilon} \id_{V\otimes V}.
			\] 
			\end{defi}
		\begin{lemm}
			The morphism $\Phi_{\varepsilon}$ from \Cref{map from hecke to endos} factors through $\TL(B_n)$ 
			\begin{equation} \label{comm TL diagram}
				\begin{tikzcd}
					\Hecke_{1,q}(B_n) \arrow[rd, "\Phi_{\varepsilon}"] \arrow[d, twoheadrightarrow, "\varphi_{\varepsilon}"'] &                     &      \\
					\TL(B_n)_{\varepsilon} \arrow[r, "\exists!\Psi_{\varepsilon}"]                            & \End_{\Vq}(V_{\varepsilon}^{\otimes n})
				\end{tikzcd},
			\end{equation}
			where $\Psi_{\varepsilon}\colon \TL(B_n)_{\varepsilon} \rightarrow \End_{\C(q)}(V_{\varepsilon}^{\otimes n}) $ is given by
			\begin{equation*}
				s_0 \mapsto \kappa_{V_{\varepsilon}} \otimes \id_{V^{\otimes (n-1)}}, \quad
				U_i \mapsto \id_{V^{\otimes (i-1)}} \otimes \underline{H}_{V_\varepsilon,V_\varepsilon} \otimes \id_{V^{\otimes (n-i-1)}} \text{ where } i=1,\ldots,n-1.
			\end{equation*}
		\end{lemm}
		\begin{theo}[Schur--Weyl duality] \label{Iso of TL and endo algebra}
			The map $
				\Psi_{\varepsilon}\colon \TL(B_n)_{\varepsilon}\rightarrow \End_{\Vq}(V_{\varepsilon}^{\otimes n})
			$ is an isomorphism of $\C(q)$-algebras.
		\end{theo}
		\begin{proof}
			First assume $\varepsilon=1$. Using Schur--Weyl duality \cite[Theorem 10.4]{stroppel2018} between $\Vq$ and $\Hecke_{1,q}(B_n)$, the morphism $\Phi_{1}$ is surjective and hence by the commutativity of  \eqref{comm TL diagram} the map $\Psi_{1}$ is also surjective. In order to show injectivity we compare dimensions. By the explicit decomposition into isotypical components in \Cref{Theorem decomposition of tensor power of V}, the dimension of $\End_{\Vq}(V_{1}^{\otimes n})$ is $\sum_{k=0}^{n}\binom{n}{k}^2$ which agrees with $\dim_k \TL(B_n)_{\varepsilon}=\binom{2n}{n}$ from \cite{green1998}. For $\varepsilon=-1$, note first that $V_{1}\cong V_{-1}\otimes k_{-1}$. Then use $V_{1}^{\otimes n}\cong V_{-1}^{\times n}\otimes k_{-1}^{\otimes n}$ and the fact that ${-} \otimes k_{-1}$ is a self-inverse functor. 
		\end{proof}
		\section{Presentations via generalized string calculus} \label{section string calculus}
		In this section we consider a monoidal category $\algreal=(\algreal,{\otimes})$ of $\Uq$-representations and a module category $\coireal=\coireal_{(\algreal, \otimes)}$ of $\Vq$-representations. We provide diagrammatic presentations, $\algdiag=(\algdiag, \otimes)$ and $\coidiag=\coidiag_{(\algdiag, \otimes)}$, for both.
		\subsection{Diagrammatics for \texorpdfstring{$\Uq$}{quantum sl(2)}} We first summarize (without proofs) easy facts.
		\begin{defi}
			Let $\algreal=\shift{V_1,k_{-1}}\subseteq \modd(\Uq)$ be the full monoidal subcategory generated by $V_1$ and $k_{-1}$.
		\end{defi}
		\begin{prop}\label{proposition rep theoretic observations}The following hold:
			\begin{enumerate} 
				\item If $L$ is an irreducible representation of $\Uq$ of type $1$, then $L\otimes k_{-1}$ is an irreducible represenation of type $-1$.
				\item Let $n\geq1$, $\varepsilon,\varepsilon'\in\{1,-1\}$.
				The Clebsch--Gordan rule holds:
				\[
					V_{n,\varepsilon}\otimes V_{\varepsilon'}\cong V_{n+1,\varepsilon\varepsilon'}\oplus V_{n-1,\varepsilon\varepsilon'}.
				\] 
				\item The following map defines an isomorphism of $\Uq$-representations:
				\begin{align*}
					V_\varepsilon\otimes k_{-1} &\to k_{-1} \otimes V_{\varepsilon}, &
					x_{\varepsilon} \otimes 1_{-1} &\mapsto 1_{-1} \otimes x_{\varepsilon},
					y_{\varepsilon} \otimes 1_{-1} \mapsto - 1_{-1} \otimes y_{\varepsilon}.
				\end{align*}
				\item Let $\varepsilon_1,\ldots, \varepsilon_r, \varepsilon_1',\ldots, \varepsilon_s'\in\{-1,1\}$. We have as vector spaces
				\[
				\Hom(V_{\varepsilon_1}\otimes \cdots \otimes V_{\varepsilon_r}, V_{\varepsilon_{1}'}\otimes \cdots \otimes V_{\varepsilon_{s}'})\!\cong\!\begin{cases}
					\Hom(V_1^{\otimes r}, V_1^{\otimes s})& \text{if } \varepsilon_1\cdots\varepsilon_{r}=\varepsilon_{1}'\cdots \varepsilon_{s}' \\
					0 & \text{else.}
				\end{cases}
				\]
			\end{enumerate}
		\end{prop}
		\begin{defi}
			Let $\algdiag$ be the $k$-linear monoidal category generated by two objects $\redbullet$ and $\blueminus$, and morphisms
			\begin{gather*}
				\cbox{
					\begin{tikzpicture}[tldiagram, scale=2/3]
						I\draw[red] \tlcoord{0}{0} \capright;
					\end{tikzpicture}
				}\colon \redbullet\redbullet\to \mU, \,	\cbox{
					\begin{tikzpicture}[tldiagram, scale=2/3]
						I\draw[red] \tlcoord{0}{0} \cupright;
					\end{tikzpicture}
				}\colon  \mU \to \redbullet\redbullet, \, \cbox{
					\begin{tikzpicture}[tldiagram, scale=2/3]
						I\draw[blue] \tlcoord{0}{0} \capright;
					\end{tikzpicture}
				}\colon \blueminus\blueminus\to \mU, \, \cbox{
					\begin{tikzpicture}[tldiagram, scale=2/3]
						I\draw[blue] \tlcoord{0}{0} \cupright;
					\end{tikzpicture}
				}\colon \mU\to\blueminus\blueminus, \\
				\cbox{
					\begin{tikzpicture}[tldiagram, scale=2/3]
						\posredbluecrossing{0}{0}
					\end{tikzpicture}
				}\colon \redbullet\blueminus\to \blueminus\redbullet, \,
				\cbox{
					\begin{tikzpicture}[tldiagram, scale=2/3]
						\posblueredcrossing{0}{0}
					\end{tikzpicture}
				}\colon \blueminus\redbullet\to\redbullet\blueminus, \, 
				\id_{\redbullet}=\cbox{
					\begin{tikzpicture}[tldiagram, scale=2/3]
						\draw[red] \tlcoord{0}{0} \lineup;
					\end{tikzpicture}
				}, \, \id_{\blueminus}=\cbox{
					\begin{tikzpicture}[tldiagram, scale=2/3]
						\draw[blue] \tlcoord{0}{0} \lineup;
					\end{tikzpicture}
				}
			\end{gather*}subject to the following relations:
			\begin{enumerate}
				\item \label{relation 1} invertibility of braiding: $\cbox{
					\begin{tikzpicture}[tldiagram, yscale=1/2, xscale=2/3]
						\posblueredcrossing{0}{0}
						\posredbluecrossing{1}{0}
					\end{tikzpicture}
				}=\cbox{
					\begin{tikzpicture}[tldiagram, yscale=1/2, xscale=2/3]
						\draw[blue] \tlcoord{0}{0} \lineup;
						\draw[red] \tlcoord{0}{1} \lineup;
					\end{tikzpicture}
				}, \quad \cbox{
					\begin{tikzpicture}[tldiagram, yscale=1/2, xscale=2/3]
						\posredbluecrossing{0}{0}
						\posblueredcrossing{1}{0}
					\end{tikzpicture}
				}=\cbox{
					\begin{tikzpicture}[tldiagram, yscale=1/2, xscale=2/3]
						\draw[red] \tlcoord{0}{0} \lineup;
						\draw[blue] \tlcoord{0}{1} \lineup;
					\end{tikzpicture}
				}$ 
				\item \label{relation 2} snake relations: $
				\cbox{
					\begin{tikzpicture}[tldiagram, yscale=1/2, xscale=2/3]
						\draw[red] \tlcoord{0}{0} \lineup \capright \cupright \lineup;
					\end{tikzpicture}
				} = \cbox{
					\begin{tikzpicture}[tldiagram, yscale=1/2, xscale=2/3]
						\draw[red] \tlcoord{0}{0}  \lineup;
					\end{tikzpicture}
				}=\cbox{
					\begin{tikzpicture}[tldiagram, yscale=1/2, xscale=2/3, xscale=-1]
						\draw[red] \tlcoord{0}{0} \lineup \capright \cupright \lineup;
					\end{tikzpicture}
				}, \quad \cbox{
					\begin{tikzpicture}[tldiagram, yscale=1/2, xscale=2/3]
						\draw[blue] \tlcoord{0}{0} \lineup \capright \cupright \lineup;
					\end{tikzpicture}
				} = \cbox{
					\begin{tikzpicture}[tldiagram, yscale=1/2, xscale=2/3]
						\draw[blue] \tlcoord{0}{0}  \lineup;
					\end{tikzpicture}
				}=\cbox{
					\begin{tikzpicture}[tldiagram, yscale=1/2, xscale=2/3, xscale=-1]
						\draw[blue] \tlcoord{0}{0} \lineup \capright \cupright \lineup;
					\end{tikzpicture}
				}$
\item \label{relation 3} bubble relations:
				$\cbox{
					\begin{tikzpicture}[tldiagram, yscale=1/2, xscale=2/3]
						\draw[red] \tlcoord{0}{0} \capright \cupleft;
					\end{tikzpicture}
				}=-(q+q^{-1}), \quad \cbox{
					\begin{tikzpicture}[tldiagram, yscale=1/2, xscale=2/3]
						\draw[blue] \tlcoord{0}{0} \capright \cupleft;
					\end{tikzpicture}
				}=1$
				\item \label{relation 4} $\blueminus\blueminus\cong \mU$-relation:
				$\cbox{
					\begin{tikzpicture}[tldiagram, yscale=1/2, xscale=2/3]
						\draw[blue] \tlcoord{0}{0} \capright;
						\draw[blue] \tlcoord{1}{0} \cupright;
					\end{tikzpicture}
				} = \cbox{
					\begin{tikzpicture}[tldiagram, yscale=1/2, xscale=2/3]
						\draw[blue] \tlcoord{0}{0} \lineup;
						\draw[blue] \tlcoord{0}{1} \lineup;
					\end{tikzpicture}
				}$
				\item \label{relation 5} Brauer-relations: $
				\cbox{
					\begin{tikzpicture}[tldiagram, yscale=1/2, xscale=2/3]
						\posblueredcrossing{0}{1}
						\draw[red] \tlcoord{0}{0} \lineup \capright;
						\draw[blue] \tlcoord{1}{2} \halflineup;
					\end{tikzpicture}
				} = - \cbox{
					\begin{tikzpicture}[tldiagram, yscale=1/2, xscale=2/3]
						\posredbluecrossing{0}{0}
						\draw[red] \tlcoord{0}{2} \lineup \capleft;
						\draw[blue] \tlcoord{1}{0} \halflineup;
					\end{tikzpicture}
				}, \quad \cbox{
					\begin{tikzpicture}[tldiagram, yscale=1/2, xscale=2/3]
						\posredbluecrossing{0}{1}
						\draw[blue] \tlcoord{0}{0} \lineup \capright;
						\draw[red] \tlcoord{1}{2} \halflineup;
					\end{tikzpicture}
				} =  \cbox{
					\begin{tikzpicture}[tldiagram, yscale=1/2, xscale=2/3]
						\posblueredcrossing{0}{0}
						\draw[blue] \tlcoord{0}{2} \lineup \capleft;
						\draw[red] \tlcoord{1}{0} \halflineup;
					\end{tikzpicture}
				}$
			\end{enumerate}
		\end{defi}
	\begin{prop} \label{proposition diagrammatic observations}
			The following statements about $\cD$ hold:
			\begin{enumerate}
				\item \label{prop item 1} The vertically reflected images of the Brauer relations hold:
				\[
				\cbox{
					\begin{tikzpicture}[tldiagram, yscale=1/2, xscale=2/3, yscale=-1]
						\posblueredcrossing{0}{1}
						\draw[red] \tlcoord{0}{0} \lineup \capright;
						\draw[blue] \tlcoord{1}{2} \halflineup;
					\end{tikzpicture}
				} = - \cbox{
					\begin{tikzpicture}[tldiagram, yscale=1/2, xscale=2/3, yscale=-1]
						\posredbluecrossing{0}{0}
						\draw[red] \tlcoord{0}{2} \lineup \capleft;
						\draw[blue] \tlcoord{1}{0} \halflineup;
					\end{tikzpicture}
				}, \quad \cbox{
					\begin{tikzpicture}[tldiagram, yscale=1/2, xscale=2/3, yscale=-1]
						\posredbluecrossing{0}{1}
						\draw[blue] \tlcoord{0}{0} \lineup \capright;
						\draw[red] \tlcoord{1}{2} \halflineup;
					\end{tikzpicture}
				} = \cbox{
					\begin{tikzpicture}[tldiagram, yscale=1/2, xscale=2/3, yscale=-1]
						\posblueredcrossing{0}{0}
						\draw[blue] \tlcoord{0}{2} \lineup \capleft;
						\draw[red] \tlcoord{1}{0} \halflineup;
					\end{tikzpicture}
				}
				\]
				\item \label{prop item 2} A compatibility between the snake and braid holds:
				\[
				\cbox{
					\begin{tikzpicture}[tldiagram, yscale=1/2, xscale=2/3]
						\draw[red] \tlcoord{-0.5}{0} \halflineup \lineup \capright \linewave{-1}{1} \cupright \lineup \halflineup; 
						\draw[white, double=blue] \tlcoord {-0.5}{1} \halflineup \linewave{1}{1} \halflineup;
					\end{tikzpicture}	
				}  = \cbox{
					\begin{tikzpicture}[tldiagram, yscale=1/2, xscale=2/3]
						\posredbluecrossing{0}{0}
					\end{tikzpicture}
				}, \quad \cbox{
					\begin{tikzpicture}[tldiagram, yscale=1/2, xscale=2/3]
						\draw[red] \tlcoord {-0.5}{1} \halflineup \linewave{1}{1} \halflineup;
						\draw[white, double=blue] \tlcoord{-0.5}{0} \halflineup \lineup \capright \linewave{-1}{1} \cupright \lineup \halflineup; 
					\end{tikzpicture}	
				}  = \cbox{
					\begin{tikzpicture}[tldiagram, yscale=1/2, xscale=2/3]
						\posblueredcrossing{0}{0}
					\end{tikzpicture}
				}
				\]
				\item \label{prop item 3}If $x,y$ are words in $\{\redbullet, \blueminus\}$ such that the partity modulo $2$ of $\#\{\blueminus\text{ in }x\}$ and $\#\{ \blueminus\text{ in }y\}$ are different, then $\Hom_{\algdiag}(x,y)=0$.
				\item \label{prop item 4}For a word $z$ in $\{\redbullet, \blueminus\}$, let $\bar{z}$ be $z$ with all $\blueminus$ removed. Set $z'\coloneqq \bar{z}$, if this number of removed $\blueminus$ was even and $z'\coloneqq \blueminus \bar{z}$ if the number of $\blueminus$ removed was odd. Then for all $x,y$ words in $\{\redbullet, \blueminus\}$ we have $\Hom_{\algdiag}(x,y)\cong\Hom_{\algdiag}(x',y')$.
				\item \label{prop item 5} For $n\in\mN$ let $\redbullet^n=\underbrace{\redbullet\cdots\redbullet}_{n \text{ times}}$. Let $r,s\in\mN$. Then $\Hom_{\algdiag}(\redbullet^r,\redbullet^s)$ (respectively $\Hom_{\algdiag}(\blueminus\redbullet^r,\blueminus\redbullet^s)$) is spanned by completely crossingless red matchings (respectively crossingless red matchings with a blue strand on the left).
			\end{enumerate}
		\end{prop}
		\begin{proof}
				For \ref{prop item 1} use first the Brauer relation \ref{relation 5} and apply on both sides the snake relation \ref{relation 2} each. \ref{prop item 2} follows by applying relation \ref{relation 5} to a subdiagram and then applying a snake relation. \ref{prop item 3} is clear by definition. To prove \ref{prop item 4} one can use the isomorphisms $\cbox{
					\begin{tikzpicture}[tldiagram, yscale=1/2, xscale=2/3]
						\posblueredcrossing{0}{0}
					\end{tikzpicture}
				}$, $\cbox{
					\begin{tikzpicture}[tldiagram, yscale=1/2, xscale=2/3]
						\posredbluecrossing{0}{0}
					\end{tikzpicture}
				}$, $\cbox{
					\begin{tikzpicture}[tldiagram, yscale=1/2, xscale=2/3]
						I\draw[blue] \tlcoord{0}{0} \capright;
					\end{tikzpicture}
				}, \cbox{
					\begin{tikzpicture}[tldiagram, yscale=1/2, xscale=2/3]
						I\draw[blue] \tlcoord{0}{0} \cupright;
					\end{tikzpicture}
				}$ to obtain $z\cong z'$ in $\algdiag$. Hence $\Hom_{\algdiag}(x,y)\cong\Hom(x',y')$. We divide the proof of \ref{prop item 5} into two claims: \\ 
				\underline{Claim 1:} We can assume without loss of generality that the diagram we consider contains no closed components $C$. To see this take a diagram $d$ which contains such a closed component $C$ (when writing $d$ as a composition of generators tensored with identities). Thus $C$ must contain cup and cap. We do a case by case distinction:
				\emph{Case 1:} $C$ is a small circle $\tlcircle$. Using bubble removal, $C$ can be replaced by a scalar.
				\emph{Case 2:} $C$ contains a kink $\cbox{
					\begin{tikzpicture}[scale=0.07]
						\draw \tlcoord{0}{0} \lineup \capright \cupright \lineup;
					\end{tikzpicture}
				}$ or $\cbox{
					\begin{tikzpicture}[scale=0.07, xscale=-1]
						\draw \tlcoord{0}{0} \lineup \capright \cupright \lineup;
					\end{tikzpicture}
				}$. Then replace $C$ by $\tlline$ using the snake relations. \emph{Case 3:} The only remaining possibility is that the closed component gets crossed by some number of strands of a different color. Using height moves, the invertibility of the braiding and the Brauer relation repeatedly we can assume the strands that cross $C$, start on one side of $C$ and end on the different side of $C$. For instance:
				\[
				\cbox
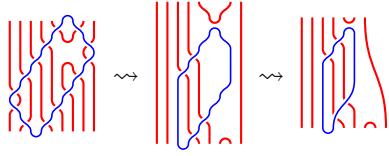
{
					\begin{tikzpicture}[tldiagram, scale=0.2]
						\draw[red] \tlcoord{0}{0} \lineup \linewave{1}{1} \linewave{1}{-1} \lineup \lineup \lineup \lineup;
						\draw[red] \tlcoord{0}{1} \linewave{1}{1} \lineup \lineup \linewave{1}{-1} \lineup \lineup \lineup;
						\draw[red] \tlcoord{0}{4} \linewave{1}{-1} \lineup \lineup \lineup \linewave{1}{-1} \lineup \lineup;
						\draw[red] \tlcoord{0}{5} \lineup \linewave{1}{-1} \lineup \lineup \lineup \linewave{1}{-1} \lineup;
						\draw[red] \tlcoord{0}{6} \lineup \lineup \linewave{1}{-1} \lineup \capright \linewave{-1}{1} \linedown \linedown \linedown;
						\draw[red] \tlcoord{0}{8} \lineup \lineup \lineup \lineup \linewave{1}{-1} \linewave{1}{1} \lineup;
						\draw[red] \tlcoord{7}{4} \linewave{-1}{1} \cupright \linewave{1}{1}; 
						\draw[white, double=blue] \tlcoord{0}{2} \linewave{1}{-1} \linewave{1}{-1} \linewave{1}{1} \linewave{1}{1}\linewave{1}{1}\linewave{1}{1}\linewave{1}{1}\capright \linewave{-1}{1} \linewave{-1}{1} \linewave{-1}{-1}\linewave{-1}{-1}\linewave{-1}{-1}\linewave{-1}{-1}\linewave{-1}{-1} \cupleft; 
					\end{tikzpicture}	
				} \rightsquigarrow \cbox{
					\begin{tikzpicture}[tldiagram, scale=0.2]
						\draw[red] \tlcoord{0}{0} \lineup \lineup \lineup \lineup \lineup \lineup \lineup \lineup \lineup;
						\draw[red] \tlcoord{0}{1} \lineup \lineup \lineup \lineup \lineup \lineup \lineup \lineup \lineup;
						\draw[red] \tlcoord{0}{4} \linewave{1}{-1} \lineup \lineup \lineup \linewave{1}{-1} \lineup \lineup \lineup \lineup;
						\draw[red] \tlcoord{0}{5} \lineup \linewave{1}{-1} \lineup \lineup \lineup \linewave{1}{-1} \lineup \lineup \lineup;
						\draw[red] \tlcoord{0}{6} \capright   ;
						\draw[red] \tlcoord{0}{8} \lineup \lineup \lineup \lineup \lineup \lineup \lineup \lineup \lineup;
						\draw[red] \tlcoord{9}{4} \linewave{-1}{1} \cupright \linewave{1}{1}; 
						\draw[white, double=blue] \tlcoord{0}{2} \lineup \lineup \lineup \lineup \linewave{1}{1}\linewave{1}{1}\linewave{1}{1}\capright \linewave{-1}{1} \linedown \linedown \linewave{-1}{-1}\linewave{-1}{-1}\linewave{-1}{-1}\linewave{-1}{-1} \cupleft; 
					\end{tikzpicture}	
				} \rightsquigarrow \cbox{
					\begin{tikzpicture}[tldiagram, scale=0.2]
						\draw[red] \tlcoord{0}{0} \lineup \lineup \lineup \lineup \lineup \lineup \lineup ;
						\draw[red] \tlcoord{0}{1} \lineup \lineup \lineup \lineup \lineup \lineup \lineup ;
						\draw[red] \tlcoord{0}{4} \linewave{1}{-1} \lineup \lineup \lineup \linewave{1}{-1}  \lineup \lineup;
						\draw[red] \tlcoord{0}{5} \lineup \linewave{1}{-1} \lineup \lineup \lineup \linewave{1}{-1} \lineup ;
						\draw[red] \tlcoord{0}{6} \capright   ;
						\draw[red] \tlcoord{0}{8} \linewave{7}{-2};
						\draw[red] \tlcoord{7}{4} \cupright; 
						\draw[white, double=blue] \tlcoord{0}{2} \lineup \lineup \lineup \lineup \linewave{1}{1}\linewave{1}{1}\capright \linedown \linedown \linedown \linewave{-3}{-2} \cupleft; 
					\end{tikzpicture}	
				}
				\]
				Now the claim follows by induction on the number of strands crossing $C$ (in the top right picture this number is two) using the Brauer relation and the invertibility of the braiding. \\
				\underline{Claim 2:} Assume $d$ is a diagram in this $\Hom$-space. Then one of the following holds.
				\begin{enumerate}[label=(\alph*)]
					\item  The diagram is entirely red, i.e.\ there is no blue component at all. (In this case the statement is just about the usual Temperley--Lieb category.) \label{item a}
					\item  We can rewrite $d$ as 
					\[
					\cbox
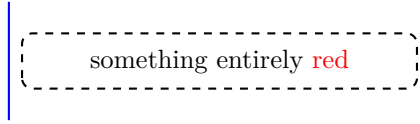
{
						\begin{tikzpicture}[tldiagram]
							\draw[blue] \tlcoord{-0.75}{0} \halflineup \lineup;
							\draw[dashed] \tlcoord{0}{1}  \maketlboxnormal{7}{something entirely \textcolor{red}{red}};
						\end{tikzpicture}	
					} \label{item b}
					\]
				\end{enumerate}
				In case $d\in \Hom_{\algdiag}(\redbullet^r,\redbullet^s)$ \ref{item a} follows by Claim 1. If $d\in \Hom_{\algdiag}(\blueminus\redbullet^r,\blueminus\redbullet^s)$ we have only one blue component by Claim 1, which is a strand $S$ connecting the two $\blueminus$. Let $m$ be the maximum of of strands which appear at some level to the left of $S$. If $m=0$, we are done. If $m>0$, look at the first time this appears (reading from top to bottom). After applying a height move we end up in one of the following two situations:
				\begin{align*}
					\cbox
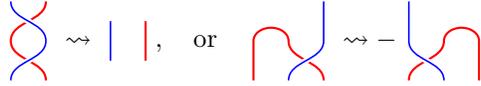
{
						\begin{tikzpicture}[tldiagram, yscale=1/2, xscale=2/3]
							\draw[red] \tlcoord{0}{1}\linewave{1}{-1} \linewave{1}{1};
							\draw[white, double=blue] \tlcoord{0}{0}\linewave{1}{1} \linewave{1}{-1};
						\end{tikzpicture}
					}  \rightsquigarrow \cbox{
						\begin{tikzpicture}[tldiagram, yscale=1/2, xscale=2/3]
							\draw[red] \tlcoord{0}{1} \lineup;
							\draw[white, double=blue] \tlcoord{0}{0} \lineup; 
						\end{tikzpicture}
					}, \quad{\text{or}} \quad\cbox{
						\begin{tikzpicture}[tldiagram, yscale=1/2, xscale=2/3]
							\draw[red] \tlcoord{0}{0}\lineup \capright \linewave{-1}{1};
							\draw[white, double=blue] \tlcoord{0}{1}\linewave{1}{1} \lineup;
						\end{tikzpicture}
					}  \rightsquigarrow - \cbox{
						\begin{tikzpicture}[tldiagram, yscale=1/2, xscale=2/3]
							\draw[red] \tlcoord{0}{0} \linewave{1}{1} \capright \linedown;
							\draw[white, double=blue] \tlcoord{0}{1} \linewave{1}{-1} \lineup;
						\end{tikzpicture}.
					}
				\end{align*}
				In any case, we reduce $m$ or at least reduce the number of times $m$ is achieved. Repeating this argument gives finally $m=0$, which shows \ref{item b}. 
				
				Now we can apply the usual type $A$ Temperley--Lieb category diagrammatics to the completely red part, and \ref{prop item 5} follows.
		\end{proof}
		\begin{defi}
			We define morphisms of $\Uq$-modules 
			\begin{gather*}
				\textcolor{red}{\cap}\colon V_1\otimes V_1  \rightarrow k_1, \,  \,\textcolor{red}{\cup}\colon k_1 \rightarrow V_1 \otimes V_1, \, \,
				\textcolor{blue}{\cap}\colon k_{-1}\otimes k_{-1} \rightarrow k_1, \, \,
				\textcolor{blue}{\cup}\colon k_1 \rightarrow k_{-1} \otimes k_{-1} \\
				\littleredblue\colon V_1\otimes k_{-1} \rightarrow k_{-1}\otimes V_1, \, \,
				\littlebluered\colon k_{-1}\otimes V_{1} \rightarrow V_1\otimes k_{-1}
				\shortintertext{by the formulas}
				\textcolor{red}{\cap}(x_{1} \otimes x_{1}) = \textcolor{red}{\cap}(y_{1} \otimes y_{1}) = 0, \,
				\textcolor{red}{\cap}(x_{1} \otimes y_{1}) = -q^{-1}1_{1}, \,
				\textcolor{red}{\cap}(y_{1} \otimes x_{1}) = 1_{1}, \\
				\textcolor{red}{\cup}(1_{1}) =  x_{1} \otimes y_{1} -  q y_{1} \otimes x_{1} , \, 
			\textcolor{blue}{\cap}(1_{-1} \otimes 1_{-1})=1_{1},   \,\textcolor{blue}{\cup}(1_{1})= 1_{-1} \otimes 1_{-1}. \\
			\littleredblue(x_{1} \otimes 1_{-1}) = 1_{-1} \otimes x_{1}, \,
			\littleredblue(y_{1} \otimes 1_{-1}) = -1_{-1} \otimes y_{1}, \\
			\littlebluered(1_{-1}\otimes x_1) = x_1 \otimes 1_{-1}, \,
			\littlebluered(1_{-1} \otimes y_{1}) = -y_{1} \otimes 1_{-1}.
			\end{gather*}
		\end{defi}
		\begin{theo}\label{equivV}
			The assignments
			\begin{align*}
				\Phi\colon \algdiag &\to \algreal \\
				{\textit{on objects:} }\quad \redbullet &\mapsto V_{1}, \, \blueminus \mapsto k_{-1}  \\
				{\textit{on morphisms:}}\quad \textcolor{red}{\cap} &\mapsto (\textcolor{red}{\cap}\colon V_1^{\otimes 2}\to k_1), \, 	\textcolor{red}{\cup} \mapsto (\textcolor{red}{\cup}\colon k_1 \to V_1^{\otimes 2}) \\
				\textcolor{blue}{\cap} &\mapsto (\textcolor{blue}{\cap}\colon k_{-1}^{\otimes 2}\to k_1),\, \textcolor{blue}{\cup} \mapsto (\textcolor{blue}{\cup}\colon k_1 \to k_{-1}^{\otimes 2})
				\\
				\littleredblue &\mapsto \littleredblue \colon V_1\otimes k_{-1}\to k_{-1}\otimes V_{1}, \,
				\littlebluered\mapsto \littlebluered \colon k_{-1}\otimes V_1\to V_{1}\otimes k_{-1}
			\end{align*}
			extend to an equivalence of $k$-linear monoidal categories $\Phi\colon \algdiag\to \algreal$.
		\end{theo}
		\begin{proof}\textit{Well-defined $k$-linear monoidal functor:} Relation \ref{relation 1} is clear by definition, relation \ref{relation 2} is a standard type $1$ relation, for blue this is clear by definition, \ref{relation 3} and \ref{relation 4} are obvious, for \ref{relation 5} this is an explicit calculation. \\
			\textit{Fullness:} Let $x,y\in \algreal$. By Proposition \ref{proposition rep theoretic observations} \ref{relation 4} we may assume  $\sgn{x}=\sgn{y}$. By the proof of \ref{relation 4}, any morphism  can be written as composition in \Cref{Bild} i).
			
		\begin{figure}[t]
			\[ i)\quad\quad
			\cbox{
				\begin{tikzpicture}[tldiagram, yscale=2/3]
					\draw \tlcoord{0}{2} \lineup \lineup \lineup \lineup;
					\draw \tlcoord{4}{1} \maketlboxnormal{3}{crossing isos};
					\draw \tlcoord{3}{0} \maketlboxnormal{5}{blue caps $\otimes$ identities};
					\draw \tlcoord{2}{1} \maketlboxnormal{3}{\textcolor{blue}{$\vert$} $\otimes f$ or $f$};
					\draw \tlcoord{1}{0} \maketlboxnormal{5}{blue cups $\otimes$ identities};
					\draw \tlcoord{0}{1} \maketlboxnormal{3}{crossing isos};
				\end{tikzpicture}	
			}.
			\quad\quad ii)\quad\quad
			\cbox{
				\begin{tikzpicture}[tldiagram, yscale=2/3]
					\draw \tlcoord{0}{2} \lineup \lineup \lineup \lineup \lineup \lineup;
					\draw \tlcoord{6}{1} \maketlboxnormal{3}{crossing isos};
					\draw \tlcoord{5}{0} \maketlboxnormal{5}{blue caps $\otimes$ identities};
					\draw \tlcoord{4}{0} \maketlboxnormal{5}{possibly $\textcolor{blue}{\perp}\otimes$ identities};
					\draw \tlcoord{3}{1.5} \maketlboxnormal{2}{$f$};
					\draw \tlcoord{2}{0} \maketlboxnormal{5}{possibly $\textcolor{blue}{\top}$ $\otimes$  identities};
					\draw \tlcoord{1}{0} \maketlboxnormal{5}{blue cups $\otimes$ identities};
					\draw \tlcoord{0}{1} \maketlboxnormal{3}{crossing isos};
				\end{tikzpicture}	
			}
			\]
			\caption{Compositions \text{with $f\in \Hom_{\Uq}(V_1^{\otimes r}, V_1^{\otimes s})$ for some $r,s$}}\label{Bild}
			\end{figure}
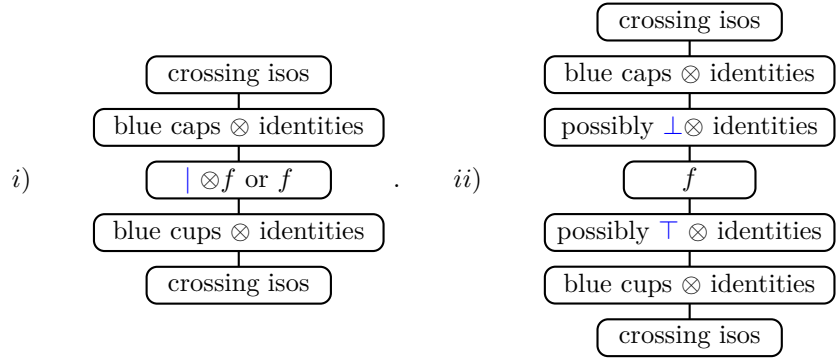
			Hence fullness follows from the analogous statement for type $1$, which gives that $f$ lies in the image of the usual Temperley--Lieb category. \\
			\textit{Faithfulness:} Since the vector space generators from \Cref{proposition diagrammatic observations}, \ref{prop item 5} are sent bijectively to the type $A$ Temperley--Lieb basis of $\End_{\Uq}(V_1^{\otimes n})$, they are linearly independent and $\Phi$ is faithful.
		\end{proof}
		
		\begin{rema} \label{remark category D is braided}
			One obtains an action of the braid group on $\algdiag$ via the morphisms
			\[
			\cbox{
				\begin{tikzpicture}[tldiagram, scale=2/3]
					\posredbluecrossing{0}{1}
				\end{tikzpicture}
			}, \quad
			\cbox{
				\begin{tikzpicture}[tldiagram, scale=2/3]
					\posblueredcrossing{0}{1}
				\end{tikzpicture}
			}, \quad
			\cbox{
				\begin{tikzpicture}[tldiagram, scale=2/3]
					\posredredcrossing{0}{1}
				\end{tikzpicture}
			}\coloneqq\cbox{
				\begin{tikzpicture}[tldiagram, scale=2/3]
					\draw[red] \tlcoord{0}{0} \capright;
					\draw[red] \tlcoord{1}{0} \cupright;
				\end{tikzpicture}
			} + q^{-1} \cbox{
				\begin{tikzpicture}[tldiagram, scale=2/3]
					\draw[red] \tlcoord{0}{0} \lineup;
					\draw[red] \tlcoord{0}{1} \lineup;
				\end{tikzpicture}
			}, \quad \cbox{
				\begin{tikzpicture}[tldiagram, scale=2/3]
					\posbluebluecrossing{0}{0}
				\end{tikzpicture}
			}\coloneqq\cbox{
				\begin{tikzpicture}[tldiagram, scale=2/3]
					\draw[blue] \tlcoord{0}{0} \lineup;
					\draw[blue] \tlcoord{0}{1} \lineup;
				\end{tikzpicture}
			}.
			\]
		\end{rema}
\subsection{Diagrammatics for \texorpdfstring{$\Vq$}{the coideal}}	
		\begin{defi}
			Let $\coireal=\shift{k_1}\algreal\subseteq \modd(\Vq)$ be the right module full subcategory generated by $k_{1}=\mU$ over the monoidal category $\algreal\subseteq \modd(\Uq)$, see \cite{HaOl-actions-tensor-categories}. By definition this means we have a commutative diagram:
			\[
			\begin{tikzcd}
				\coireal \times \algreal \arrow[d, "\text{inclusion}"', hook] \arrow[r, "\text{action}"] & \coireal \arrow[d, "\text{inclusion}", hook] \\
				\modd(\Vq)\times \modd(\Uq) \arrow[r] \arrow[r, "\text{action}"]                & \modd(\Vq)                              
			\end{tikzcd}
			\]
		\end{defi}
		\begin{defi}
			Let $(\coidiag, \otimes\colon \coidiag\times \algdiag \to \coidiag)$ be the $k$-linear right module category over $\algdiag$ generated by one object $\mU'$ and morphisms
			\[
			\cbox{
				\begin{tikzpicture}[tldiagram, scale=2/3]
					\draw[dotted] \tlcoord{0}{-1} \lineup; 
					\draw[red] \tlcoord{0}{0} \dlineup; 
				\end{tikzpicture}	
			}\colon \mU' \redbullet \to \mU' \redbullet, \quad \cbox{
				\begin{tikzpicture}[tldiagram, scale=2/3]
					\draw[dotted] \tlcoord{0}{-1} \lineup; 
					\draw[blue] \tlcoord{0}{0} \linewave{0.75}{-1} \widerhakenup; 
				\end{tikzpicture}	
			}\colon \mU' \blueminus \to \mU', \quad \cbox{
				\begin{tikzpicture}[tldiagram, scale=2/3, yscale=-1]
					\draw[dotted] \tlcoord{0}{-1} \lineup; 
					\draw[blue] \tlcoord{0}{0} \linewave{0.75}{-1} \widerhakenup; 
				\end{tikzpicture}	
			}\colon    \mU'\to \mU'{\blueminus} , \quad  \id_{\mU'} = \cbox{
			\begin{tikzpicture}[tldiagram, scale=2/3]
				\draw[dotted] \tlcoord{0}{-1} \lineup;  
			\end{tikzpicture}
		}
			\]
			subject to the relations
			\begin{enumerate}[itemsep=2pt]
				\item quadratic relations: $
				\cbox{
					\begin{tikzpicture}[tldiagram, scale=2/3]
						\draw[dotted] \tlcoord{0}{-1} \lineup;
						\draw[red] \tlcoord{0.5}{0} \onedot \capright \cupleft;;
					\end{tikzpicture}	
				} = 0
				, \quad \cbox{
					\begin{tikzpicture}[tldiagram, scale=2/3]
						\draw[dotted] \tlcoord{0}{-1} \lineup; 
						\draw[red] \tlcoord{0}{0} \quarterlineup \onedot \halflineup \onedot \quarterlineup; 
					\end{tikzpicture}
				} = \cbox{
					\begin{tikzpicture}[tldiagram, scale=2/3]
						\draw[dotted] \tlcoord{0}{-1} \lineup; 
						\draw[red] \tlcoord{0}{0} \lineup; 
					\end{tikzpicture}
				}$
				\item left unit and counit: $\cbox{
					\begin{tikzpicture}[tldiagram, scale=2/3]
						\draw[dotted] \tlcoord{0}{-1} \lineup \halflineup; 
						\draw[blue] \tlcoord{0}{1} \lineup \capleft \linewave{-0.75}{-1} \widerhakendown; 
					\end{tikzpicture}
				} = \cbox{
					\begin{tikzpicture}[tldiagram, scale=2/3]
						\draw[dotted] \tlcoord{0}{-1} \lineup; 
						\draw[blue] \tlcoord{0}{0} \linewave{0.75}{-1} \widerhakenup; 
					\end{tikzpicture}	
				}, \quad \cbox{
					\begin{tikzpicture}[tldiagram, scale=2/3, yscale=-1]
						\draw[dotted] \tlcoord{0}{-1} \lineup \halflineup; 
						\draw[blue] \tlcoord{0}{1} \lineup \capleft \linewave{-0.75}{-1} \widerhakendown; 
					\end{tikzpicture}
				} = \cbox{
					\begin{tikzpicture}[tldiagram, scale=2/3, yscale=-1]
						\draw[dotted] \tlcoord{0}{-1} \lineup; 
						\draw[blue] \tlcoord{0}{0} \linewave{0.75}{-1} \widerhakenup; 
					\end{tikzpicture}	
				}$
				\item $\mU'\cong \mU'\otimes \blueminus$-relations: $\cbox{
					\begin{tikzpicture}[tldiagram, scale=2/3]
						\draw[dotted] \tlcoord{0}{-1} \lineup \lineup; 
						\draw[blue] \tlcoord{0}{0} \linewave{0.75}{-1} \widerhakenup; 
						\draw[blue] \tlcoord{1}{-1} \widerhakendown \linewave{0.75}{1}; 
					\end{tikzpicture}
				} = \cbox{
					\begin{tikzpicture}[tldiagram, scale=2/3]
						\draw[dotted] \tlcoord{0}{-1} \lineup; 
						\draw[blue] \tlcoord{0}{0} \lineup;
					\end{tikzpicture}
				}, \quad \cbox{
					\begin{tikzpicture}[tldiagram, scale=2/3]
						\draw[dotted] \tlcoord{0}{-1} \lineup \lineup; 
						\draw[blue] \tlcoord{0.25}{-1} \widerhakendown  \linewave{0.75}{1} \linewave{0.75}{-1} \widerhakenup; 
					\end{tikzpicture}
				} = \cbox{
					\begin{tikzpicture}[tldiagram, scale=2/3]
						\draw[dotted] \tlcoord{0}{-1} \lineup;  
					\end{tikzpicture}
				}$
				\item \label{wrecking ball} "wrecking ball"-relation: $\cbox{ 
					\begin{tikzpicture}[tldiagram, yscale=1/2, xscale=2/3]
						\draw[dotted] \tlcoord{0}{-1} \lineup \halflineup \lineup;  
						\draw[red] \tlcoord{1}{0} \halfdlineup;
						\draw[blue] \tlcoord{1}{1} \halflineup;
						\posblueredcrossing{0}{0}
						\posredbluecrossing{1.5}{0}
					\end{tikzpicture}
				} = - \cbox{
					\begin{tikzpicture}[tldiagram, yscale=1/2, xscale=2/3]
						\draw[dotted] \tlcoord{0}{-1} \lineup \lineup;  
						\draw[blue] \tlcoord{0}{0} \linewave{0.75}{-1} \widerhakenup; 
						\draw[blue] \tlcoord{1.25}{-1} \widerhakendown \linewave{0.75}{1}; 
						\draw[red] \tlcoord{0}{1} \linewave{0.75}{-1} \halfdlineup \linewave{0.75}{1};
					\end{tikzpicture}
				}$
			\end{enumerate}
		\end{defi}
\begin{rema}
			The specialized quadratic relations imply the reflection equation:
			\[
			\cbox{
				\begin{tikzpicture}[tldiagram, yscale=1/2, xscale=2/3]
					\draw[dotted] \tlcoord{0}{-1} \lineup \lineup \lineup; 
					\draw[red] \tlcoord{0}{0} \halfdlineup;
					\draw[red] \tlcoord{1.5}{0} \quarterlineup \onedot \quarterlineup;
					\draw[red] \tlcoord{0}{1} \halflineup;
					\draw[red] \tlcoord{1.5}{1} \halflineup;
					\posredredcrossing{0.5}{0}
					\posredredcrossing{2}{0}
				\end{tikzpicture}
			} = \cbox{
				\begin{tikzpicture}[tldiagram, scale=2/3]
					\draw[dotted] \tlcoord{0}{-1} \lineup; 
					\draw[red] \tlcoord{0}{0} \lineup;
					\draw[red] \tlcoord{0}{1} \lineup;
				\end{tikzpicture}
			} + q^{-1} \cbox{
				\begin{tikzpicture}[tldiagram, scale=2/3]
					\draw[dotted] \tlcoord{0}{-1} \lineup; 
					\draw[red] \tlcoord{0}{0} \capright;
					\draw[red] \tlcoord{1}{0} \cupright;
				\end{tikzpicture}
			} + q^{-1} \cbox{
				\begin{tikzpicture}[tldiagram, scale=2/3]
					\draw[dotted] \tlcoord{0}{-1} \lineup; 
					\draw[red] \tlcoord{0}{0} \dcapright;
					\draw[red] \tlcoord{1}{0} \dcupright;
				\end{tikzpicture}
			} =
			\cbox{
				\begin{tikzpicture}[tldiagram, yscale=1/2, xscale=2/3]
					\draw[dotted] \tlcoord{0}{-1} \lineup \lineup \lineup; 
					\draw[red] \tlcoord{1}{0} \halfdlineup;
					\draw[red] \tlcoord{2.5}{0} \halfdlineup;
					\draw[red] \tlcoord{1}{1} \halflineup;
					\draw[red] \tlcoord{2.5}{1} \halflineup;
					\posredredcrossing{0}{0}
					\posredredcrossing{1.5}{0}
				\end{tikzpicture}
			}.
			\]
			This special dotted braid describing the reflection equation is the type $D_2$ full-twist. Similarly to type $A$ \cite[Thm. 2.2]{rozansky2010} its powers in the $q^{-1}$-adic norm converge towards the type $D_2$ Jones--Wenzl projector $d_{2,1}$, see Theorem \ref{type D projector recursive description} for the explicit formula. For a precise statement and proof of the convergence see \cite{Woj2022}.
		\end{rema}
		\begin{theo} \label{final theorem}
			There is a $k$-linear equivalence of categories $\Phi'\colon \coidiag\to \coireal$, which is a bijection on objects such that the following diagram commutes
			\[
			\begin{tikzcd}
				\coidiag\times \algdiag \arrow[r, "\text{action}"] \arrow[d, "\Phi'\times \Phi"'] & \coidiag \arrow[d, "\Phi'"] \\
				\coireal\times \algreal \arrow[r, "\text{action}"]                                & \coireal.                   
			\end{tikzcd}
			\]
		It defines an equivalence of module categories via the equivalence from \Cref{equivV}.
		\end{theo}
		\begin{proof}[Proof of \Cref{final theorem}]
			We define $\Phi'\colon \coidiag \to \coireal$ by mapping $\mU'\mapsto k_{1}$ and
			\[
			\cbox{
				\begin{tikzpicture}[tldiagram, scale=2/3]
					\draw[dotted] \tlcoord{0}{-1} \lineup; 
					\draw[red] \tlcoord{0}{0} \dlineup; 
				\end{tikzpicture}	
			}\mapsto \id_{k_1}\otimes \kappa_{V_{1}}, \,  \, \cbox{
				\begin{tikzpicture}[tldiagram, scale=2/3]
					\draw[dotted] \tlcoord{0}{-1} \lineup; 
					\draw[blue] \tlcoord{0}{0} \linewave{0.75}{-1} \widerhakenup; 
				\end{tikzpicture}	
			} \mapsto [1_1\otimes 1_{-1}\mapsto 1_1], \, \,  \cbox{
				\begin{tikzpicture}[tldiagram, scale=2/3, yscale=-1]
					\draw[dotted] \tlcoord{0}{-1} \lineup; 
					\draw[blue] \tlcoord{0}{0} \linewave{0.75}{-1} \widerhakenup; 
				\end{tikzpicture}	
			} \mapsto [1_1\mapsto 1_1\otimes 1_{-1}].
			\]
			\textit{Functor of module categories:} All the images define morphisms of $\Vq$-modules and checking the relations is an easy calculation. \\
			
			\textit{Fullness:} Let $x,y\in \coireal$, i.e. tensor products of $V_{1}$ and $k_{-1}$. We use the braiding from \Cref{remark category D is braided} and blue cup and cap to replace $x, y$ with $\Uq$-isomorphic objects $x'\cong x$, $y'\cong y$ with at most one tensor factor $k_{-1}$ on the very left. Concretely $x'=V_1^{\otimes r}$ (or $x'=k_{-1}\otimes V_1^{\otimes r}$) and $y'=V_1^{\otimes s}$ (or $y'=k_{-1}\otimes V_1^{\otimes s}$) for some $r,s\geq 0$. Then use the $\Vq$-equivariant isomorphisms $\bluewiderhakenup,\bluewiderhakendown$ to replace, if necessary, the left tensor factor $k_{-1}$ with $k_1$. This allows to write any $\Vq$-equivariant $f\colon x\to y$ as composition as in \Cref{Bild} ii).
			The space $\Hom_{\Vq}(V_1^{\otimes r}, V_1^{\otimes s})$ becomes zero, if $r+s$ is odd, by Theorem \ref{Theorem decomposition of tensor power of V}. Hence we can assume $r+s=2n$ for some $n\geq 0$. After using that $V_1$ is self-dual via the cup and cap morphisms we get from an adjunction
			\[
			\Hom_{\Vq}(V_1^{\otimes r}, V_1^{\otimes s})\cong \Hom_{\Vq}(V_1^{\otimes n}, V_1^{\otimes n})\cong \TL(B_n)_1.
			\]
			Since all the generators of $\TL(B_n)_1$ lie in the image of the functor and all the isomorphisms we used are in the image, we obtain fullness. \\
			\textit{Faithfulness:} First notice that every object in $\coidiag$ is isomorphic to some tensor power of $\redbullet$ by diagrammatic analogues of the isomorphisms as above used for $\coireal$. Hence it is enogh to show for every $r,s\geq 0$ the faithfullness of $\Phi'$ on the $\Hom$-space $\Hom_{\coidiag}(\redbullet^r,\redbullet^s)$.
			If $r+s$ is odd, we have by choice of generators $\Hom_{\coidiag}(\redbullet^r,\redbullet^s)=0$, so faithfulness is clear.
			If $r+s=2n$ is even, use the self-duality of $\redbullet$ to replace $\Hom_{\coidiag}(\redbullet^r,\redbullet^s)$ with $\Hom_{\coidiag}(\redbullet^n,\redbullet^n)$.  
			In order to conclude the theorem one shows analogously to the proof of \Cref{proposition diagrammatic observations} v) that that every blue closed component involving $\bluewiderhakenup$ and $\bluewiderhakendown$ can be removed, the crucial relation being the wrecking ball \ref{wrecking ball}. Hence $\Hom_{\coidiag}(\redbullet^n,\redbullet^n)$ is spanned by completely red type $B$ Temperley--Lieb diagrams, which form a basis of the image under the functor by \Cref{Iso of TL and endo algebra}. 
		\end{proof}

		\section{Different types of Jones--Wenzl projectors} \label{section on type B projectors}

		 We next connect the direct sum decomposition from \Cref{section on fusion rules} with the diagrammatics from \Cref{section string calculus}. This means we study the Karoubian envelope of $\coidiag\simeq \coireal$, or more elementary speaking, idempotents in $\TL(B_n)_{1}= \End_{\coidiag}(\redbullet^n)\cong \End_{\coireal}(V_1^{\otimes n})$. Let $\varepsilon\in\{1,-1\}, n\geq 1$. Recall the Jones--Wenzl projectors, see e.g. \cite{flath95}[\S 3.5]:
		%
		
		\begin{defi} \label{recursive formula type A}
			The type $A_{n-1}$ Jones--Wenzl projector $a_{n-1,\varepsilon}$ is visualized as an uncolored (white) box on $n$ strands and determined by the recursive description
			\begin{gather*}
				a_{n-1,\varepsilon}=\cbox{
					\begin{tikzpicture}[tldiagram, yscale=1/2, xscale=2/3]
						\draw \tlcoord{-1.5}{0} \lineup \lineup \lineup;
						\draw \tlcoord{-1.5}{1} \lineup  \lineup \lineup;
						\draw \tlcoord{-1.5}{2} \lineup \lineup \lineup;
						\draw \tlcoord{0}{0} \maketlboxnormal{3}{$n, \varepsilon$};
						=
					\end{tikzpicture}
				} , \,
				\cbox{\begin{tikzpicture}[tldiagram, yscale=1/2, xscale=2/3]
						\draw \tlcoord{0}{0} \lineup \lineup;
						\draw \tlcoord{1}{0} \maketlboxnormal{1}{$1, \varepsilon$};
				\end{tikzpicture}}=\cbox{
					\begin{tikzpicture}[tldiagram]
						\draw \tlcoord{0}{0} \lineup;
				\end{tikzpicture}} \, , \, \cbox{
					\begin{tikzpicture}[tldiagram, yscale=1/2, xscale=2/3]
						\draw \tlcoord{-1.5}{0} \lineup \lineup \lineup;
						\draw \tlcoord{-1.5}{1} \lineup \lineup \lineup;
						\draw \tlcoord{-1.5}{2} \lineup \lineup \lineup;
						\draw \tlcoord{0}{0} \maketlboxnormal{3}{$n+1, \varepsilon$};
						=
					\end{tikzpicture}
				}
				=
				\cbox{
					\begin{tikzpicture}[tldiagram, yscale=1/2, xscale=2/3]
						\draw \tlcoord{-1.5}{0} \lineup \lineup \lineup;
						\draw \tlcoord{-1.5}{1} \lineup \lineup \lineup;
						\draw \tlcoord{-1.5}{2} \lineup  \lineup \lineup;
						\draw \tlcoord{0}{0} \maketlboxnormal{2}{$n, \varepsilon$};
						=
					\end{tikzpicture}
				}
				+
				\varepsilon \frac{[n]}{[n+1]}
				\cbox{
					\begin{tikzpicture}[tldiagram, yscale=1/2, xscale=2/3]
						\draw \tlcoord{-1}{0} \lineup \lineup \halflineup \lineup \lineup;
						\draw \tlcoord{2.5}{1} \lineup;
						\draw \tlcoord{0}{1} \linedown;
						\draw \tlcoord{0}{1} \smalllineup \capright \smalllinedown \linedown;
						\draw \tlcoord{2.5}{1} \smalllinedown \cupright \smalllineup \lineup;
						\draw \tlcoord{0}{0} \maketlboxnormal{2}{$n, \varepsilon$};
						\draw \tlcoord{2.5}{0} \maketlboxnormal{2}{$n, \varepsilon$};
					\end{tikzpicture}
				}.
			\end{gather*}
		\end{defi}
		\begin{defi} \label{defi type B projectors} \label{definition of projectors of type d} 
			We define elements $b_{n,{+}, \varepsilon}, b_{n,{-}, \varepsilon}\in \TL(B_n)_{\varepsilon}$ recursively and visualize them as red respectively\ blue boxes on $n$ strands:
			\[
			b_{n,{+},\varepsilon}=\cbox{
				\begin{tikzpicture}[tldiagram, yscale=1/2, xscale=2/3]
					\draw \tlcoord{-1.5}{0} \lineup \lineup \lineup;
					\draw \tlcoord{-1.5}{1} \lineup  \lineup \lineup;
					\draw \tlcoord{-1.5}{2} \lineup \lineup \lineup;
					\draw \tlcoord{0}{0} \maketlboxred{3}{$n, \varepsilon$};
					=
				\end{tikzpicture}
			}, \quad b_{n,{-}, \varepsilon} = \cbox{
				\begin{tikzpicture}[tldiagram, yscale=1/2, xscale=2/3]
					\draw \tlcoord{-1.5}{0} \lineup \lineup \lineup;
					\draw \tlcoord{-1.5}{1} \lineup \lineup  \lineup;
					\draw \tlcoord{-1.5}{2} \lineup \lineup \lineup;
					\draw \tlcoord{0}{0} \maketlboxblue{3}{$n, \varepsilon$};
					=
				\end{tikzpicture}
			}.
			\]
			First we define for $n=1$
			\[
			\cbox{\begin{tikzpicture}[tldiagram, yscale=1/2, xscale=2/3]
					\draw \tlcoord{0}{0} \lineup \lineup;
					\draw \tlcoord{1}{0} \maketlboxred{1}{$1, \varepsilon$};
			\end{tikzpicture}}\coloneqq\frac{1}{2}(1+s_0) = \frac{1}{2} \left( \,
			\cbox{\begin{tikzpicture}[tldiagram, scale=2/3]
					\draw \tlcoord{0}{0} \lineup;
			\end{tikzpicture}}
			 +  \cbox{\begin{tikzpicture}[tldiagram, scale=2/3]
					\draw \tlcoord{0}{0} \dlineup;
			\end{tikzpicture}} \, \right) \, , 	\quad  \cbox{\begin{tikzpicture}[tldiagram, yscale=1/2, xscale=2/3]
					\draw \tlcoord{0}{0} \lineup \lineup;
					\draw \tlcoord{1}{0} \maketlboxblue{1}{$1, \varepsilon$};
			\end{tikzpicture}}\coloneqq\frac{1}{2}(1-s_0) = \frac{1}{2} \left( \,
			\cbox{\begin{tikzpicture}[tldiagram, yscale=1/2, xscale=2/3]
					\draw \tlcoord{0}{0} \lineup;
			\end{tikzpicture}}
			-  \cbox{\begin{tikzpicture}[tldiagram, scale=2/3]
					\draw \tlcoord{0}{0} \dlineup;
			\end{tikzpicture}} \, \right). 
			\]
			Assuming that $b_{n, {+}, \varepsilon}$ and $b_{n, {-}, \varepsilon}$ are already defined, set
			\begin{align*}
			\cbox{
				\begin{tikzpicture}[tldiagram, yscale=1/2, xscale=2/3]
					\draw \tlcoord{-1.5}{0} \lineup \lineup \lineup;
					\draw \tlcoord{-1.5}{1} \lineup \lineup \lineup;
					\draw \tlcoord{-1.5}{2} \lineup \lineup \lineup;
					\draw \tlcoord{0}{0} \maketlboxred{3}{$n+1, \varepsilon$};
					=
				\end{tikzpicture}
			}
			\,&\coloneqq\,
			\cbox{
				\begin{tikzpicture}[tldiagram, yscale=1/2, xscale=2/3]
					\draw \tlcoord{-1.5}{0} \lineup \lineup \lineup;
					\draw \tlcoord{-1.5}{1} \lineup \lineup \lineup;
					\draw \tlcoord{-1.5}{2} \lineup  \lineup \lineup;
					\draw \tlcoord{0}{0} \maketlboxred{2}{$n, \varepsilon$};
					=
				\end{tikzpicture}
			}
			\,+\,
			\varepsilon \frac{q^{n-1} + q^{{-}(n-1)}}{q^n + q^{{-}n}}
			\cbox{
				\begin{tikzpicture}[tldiagram, yscale=1/2, xscale=2/3]
					\draw \tlcoord{-1}{0} \lineup \lineup \halflineup \lineup \lineup;
					\draw \tlcoord{2.5}{1} \lineup;
					\draw \tlcoord{0}{1} \linedown;
					\draw \tlcoord{0}{1} \smalllineup \capright \smalllinedown \linedown;
					\draw \tlcoord{2.5}{1} \smalllinedown \cupright \smalllineup \lineup;
					\draw \tlcoord{0}{0} \maketlboxred{2}{$n, \varepsilon$};
					\draw \tlcoord{2.5}{0} \maketlboxred{2}{$n, \varepsilon$};
				\end{tikzpicture}
			}, \\
		\cbox{
			\begin{tikzpicture}[tldiagram, yscale=1/2, xscale=2/3]
				\draw \tlcoord{-1.5}{0} \lineup \lineup \lineup;
				\draw \tlcoord{-1.5}{1} \lineup \lineup \lineup;
				\draw \tlcoord{-1.5}{2} \lineup \lineup \lineup;
				\draw \tlcoord{0}{0} \maketlboxblue{3}{$n+1, \varepsilon$};
				=
			\end{tikzpicture}
		}
		\,&\coloneqq\,
		\cbox{
			\begin{tikzpicture}[tldiagram, yscale=1/2, xscale=2/3]
				\draw \tlcoord{-1.5}{0} \lineup \lineup \lineup;
				\draw \tlcoord{-1.5}{1} \lineup \lineup \lineup;
				\draw \tlcoord{-1.5}{2} \lineup \lineup \lineup;
				\draw \tlcoord{0}{0} \maketlboxblue{2}{$n, \varepsilon$};
				=
			\end{tikzpicture}
		}
		\,+\,
		\varepsilon \frac{q^{n-1} + q^{{-}(n-1)}}{q^n + q^{{-}n}}
		\cbox{
			\begin{tikzpicture}[tldiagram, yscale=1/2, xscale=2/3]
				\draw \tlcoord{-1}{0} \lineup \lineup \halflineup \lineup \lineup;
				\draw \tlcoord{2.5}{1} \lineup;
				\draw \tlcoord{0}{1} \linedown;
				\draw \tlcoord{0}{1} \smalllineup \capright \smalllinedown \linedown;
				\draw \tlcoord{2.5}{1} \smalllinedown \cupright \smalllineup \lineup;
				\draw \tlcoord{0}{0} \maketlboxblue{2}{$n, \varepsilon$};
				\draw \tlcoord{2.5}{0} \maketlboxblue{2}{$n, \varepsilon$};
			\end{tikzpicture}
		}.
	\end{align*}
			We call $b_{n,{+},\varepsilon}$ (and \ $b_{n,{-},\varepsilon}$) the \textit{dominant} (respectively\ \textit{antidominant}) \textit{type $B_n$ Jones--Wenzl projector of type $\varepsilon$}. Additionally we define the \textit{type $D_n$ Jones--Wenzl projector $d_n$ of type $\varepsilon$} as the sum
			\begin{equation*}
				\cbox{
					\begin{tikzpicture}[tldiagram, yscale=1/2, xscale=2/3]
						\draw \tlcoord{-1.5}{0} \lineup \lineup  \lineup;
						\draw \tlcoord{-1.5}{1} \lineup \lineup  \lineup;
						\draw \tlcoord{-1.5}{2} \lineup \lineup  \lineup;
						\draw \tlcoord{0}{0} \maketlboxgreen{3}{$n, \varepsilon$};
						=
					\end{tikzpicture}
				} \, \coloneqq \, \cbox{
					\begin{tikzpicture}[tldiagram, yscale=1/2, xscale=2/3]
						\draw \tlcoord{-1.5}{0} \lineup \lineup \lineup;
						\draw \tlcoord{-1.5}{1} \lineup  \lineup \lineup;
						\draw \tlcoord{-1.5}{2} \lineup \lineup \lineup;
						\draw \tlcoord{0}{0} \maketlboxred{3}{$n,  \varepsilon$};
						=
					\end{tikzpicture}
				} \, + \, \cbox{
					\begin{tikzpicture}[tldiagram, yscale=1/2, xscale=2/3]
						\draw \tlcoord{-1.5}{0} \lineup \lineup \lineup;
						\draw \tlcoord{-1.5}{1} \lineup \lineup  \lineup;
						\draw \tlcoord{-1.5}{2} \lineup \lineup \lineup;
						\draw \tlcoord{0}{0} \maketlboxblue{3}{$n,  \varepsilon$};
						=
					\end{tikzpicture}
				}.
			\end{equation*}
		\end{defi}
		The red and blue projectors appeared abstractly in {\cite[\S3]{tomdieck1998}}, in their notation we have $b_{n,{+},\varepsilon}=g_n$ and $b_{n,{-},\varepsilon}=f_n$, however no diagrammatics was involved and the type $D$ projectors not defined. Our 'new' diagrammatical depiction is an attempt at generalizing the type $A$ Jones--Wenzl diagrammatics. The green projectors appeared also in non-quantized form in \cite{wedrich2018}, where quotients of affine Temperley--Lieb algebras were studied, in a context at first glance independent of type $B/D$.
		\begin{rema}
			The involution $\sigma$ from Lemma \ref{Type B automorphism s_o flips sign} swaps $b_{1,{+},\varepsilon}$ and $b_{1,{-},\varepsilon}$, thus $\sigma(b_{n,{+},\varepsilon})=b_{n,{-},\varepsilon}$ for all $n$. As a consequence, the type $D_n$ Jones--Wenzl projector $d_{n, \varepsilon}$ is contained in the fixed-point subalgebra $\TL(D_n)_{\varepsilon}^\sigma$, and thus indeed of type $D$.
		\end{rema}
		\begin{beis} \label{n=2 type B projectors}
			The recursion gives the dominant Jones--Wenzl  projector of type $B_2$
			\begin{align*}
				\cbox{
					\begin{tikzpicture}[tldiagram, yscale=1/2, xscale=2/3]
						\draw \tlcoord{0}{0} \lineup \lineup;
						\draw \tlcoord{0}{1} \lineup \lineup;
						\draw \tlcoord{1}{0} \maketlboxred{2}{$2, \varepsilon$};
					\end{tikzpicture}	
				}
				&=\frac{1}{2}\left( \cbox{
					\begin{tikzpicture}[tldiagram, scale=2/3]
						\draw \tlcoord{0}{0} \lineup;
						\draw \tlcoord{0}{1} \lineup;
					\end{tikzpicture}
				} + \cbox{
					\begin{tikzpicture}[tldiagram, scale=2/3]
						\draw \tlcoord{0}{0} \dlineup;
						\draw \tlcoord{0}{1} \lineup;
					\end{tikzpicture}
				}\right) +	\frac{\varepsilon}{2[2]} \left( 
				\cbox{
					\begin{tikzpicture}[tldiagram, scale=2/3]
						\draw \tlcoord{0}{0} \capright;
						\draw \tlcoord{1}{0} \cupright;
					\end{tikzpicture}
				} + \cbox{
					\begin{tikzpicture}[tldiagram, scale=2/3]
						\draw \tlcoord{0}{0} \capright;
						\draw \tlcoord{1}{0} \dcupright;
					\end{tikzpicture}
				}
				+ \cbox{
					\begin{tikzpicture}[tldiagram, scale=2/3]
						\draw \tlcoord{0}{0} \dcapright;
						\draw \tlcoord{1}{0} \cupright;
					\end{tikzpicture}
				}
				+ \cbox{
					\begin{tikzpicture}[tldiagram, scale=2/3]
						\draw \tlcoord{0}{0} \dcapright;
						\draw \tlcoord{1}{0} \dcupright;
					\end{tikzpicture}
				}
				\right) ,\\
			\shortintertext{and the antidominant Jones--Wenzl projector of type $B_2$}
				\cbox{
					\begin{tikzpicture}[tldiagram, yscale=1/2, xscale=2/3]
						\draw \tlcoord{0}{0} \lineup \lineup;
						\draw \tlcoord{0}{1} \lineup \lineup;
						\draw \tlcoord{1}{0} \maketlboxblue{2}{$2, \varepsilon$};
					\end{tikzpicture}	
				}&=\frac{1}{2}\left( \cbox{
					\begin{tikzpicture}[tldiagram, scale=2/3]
						\draw \tlcoord{0}{0} \lineup;
						\draw \tlcoord{0}{1} \lineup;
					\end{tikzpicture}
				} - \cbox{
					\begin{tikzpicture}[tldiagram, scale=2/3]
						\draw \tlcoord{0}{0} \dlineup;
						\draw \tlcoord{0}{1} \lineup;
					\end{tikzpicture}
				}\right) +	\frac{\varepsilon}{2[2]} \left( 
				\cbox{
					\begin{tikzpicture}[tldiagram, scale=2/3]
						\draw \tlcoord{0}{0} \capright;
						\draw \tlcoord{1}{0} \cupright;
					\end{tikzpicture}
				} - \cbox{
					\begin{tikzpicture}[tldiagram, scale=2/3]
						\draw \tlcoord{0}{0} \capright;
						\draw \tlcoord{1}{0} \dcupright;
					\end{tikzpicture}
				}
				- \cbox{
					\begin{tikzpicture}[tldiagram, scale=2/3]
						\draw \tlcoord{0}{0} \dcapright;
						\draw \tlcoord{1}{0} \cupright;
					\end{tikzpicture}
				}
				+ \cbox{
					\begin{tikzpicture}[tldiagram, scale=2/3]
						\draw \tlcoord{0}{0} \dcapright;
						\draw \tlcoord{1}{0} \dcupright;
					\end{tikzpicture}
				}
				\right).
			\end{align*}
		\end{beis}
		We will prove recursion formulas for type $D$ independent of the type $B$ result:
		\begin{theo}[Type $D$ recursion] \label{type D projector recursive description}
			The type $D$ projectors $d_{n,\varepsilon}$ are described by the following recursive formulas:
			\begin{align*}
				\cbox{\begin{tikzpicture}[tldiagram, yscale=1/2, xscale=2/3]
						\draw \tlcoord{0}{0} \lineup \lineup;
						\draw \tlcoord{1}{0} \maketlboxgreen{1}{$1, \varepsilon$};
				\end{tikzpicture}}&=\cbox{
					\begin{tikzpicture}[tldiagram, scale=2/3]
						\draw \tlcoord{0}{0} \lineup;
				\end{tikzpicture}} \, , \qquad \cbox{
					\begin{tikzpicture}[tldiagram, yscale=1/2, xscale=2/3]
						\draw \tlcoord{0}{0} \lineup \lineup;
						\draw \tlcoord{0}{1} \lineup \lineup;
						\draw \tlcoord{1}{0} \maketlboxgreen{2}{$2, \varepsilon$};
					\end{tikzpicture}	
				}=\cbox{
					\begin{tikzpicture}[tldiagram, scale=2/3]
						\draw \tlcoord{0}{0} \lineup;
						\draw \tlcoord{0}{1} \lineup;
					\end{tikzpicture}
				} + \frac{\varepsilon}{[2]} \cbox{
					\begin{tikzpicture}[tldiagram, scale=2/3]
						\draw \tlcoord{0}{0} \capright;
						\draw \tlcoord{1}{0} \cupright;
					\end{tikzpicture}
				} + \frac{\varepsilon}{[2]} \cbox{
					\begin{tikzpicture}[tldiagram, scale=2/3]
						\draw \tlcoord{0}{0} \dcapright;
						\draw \tlcoord{1}{0} \dcupright;
					\end{tikzpicture}
				} \, , \\
			\cbox{
				\begin{tikzpicture}[tldiagram, yscale=1/2, xscale=2/3]
					\draw \tlcoord{-1.5}{0} \lineup \lineup \lineup;
					\draw \tlcoord{-1.5}{1} \lineup \lineup \lineup;
					\draw \tlcoord{-1.5}{2}\lineup \lineup \lineup;
					\draw \tlcoord{0}{0} \maketlboxgreen{3}{$n+1, \varepsilon$};
					=
				\end{tikzpicture}
			}
			\,&=\,
			\cbox{
				\begin{tikzpicture}[tldiagram, yscale=1/2, xscale=2/3]
					\draw \tlcoord{-1.5}{0} \lineup \lineup \lineup;
					\draw \tlcoord{-1.5}{1} \lineup \lineup \lineup;
					\draw \tlcoord{-1.5}{2} \lineup \lineup \lineup;
					\draw \tlcoord{0}{0} \maketlboxgreen{2}{$n, \varepsilon$};
					=
				\end{tikzpicture}
			}
			\,+\,
			\varepsilon\frac{q^{n-1} + q^{{-}(n-1)}}{q^n + q^{{-}n}}
			\cbox{
				\begin{tikzpicture}[tldiagram, yscale=1/2, xscale=2/3]
					\draw \tlcoord{-1}{0} \lineup \lineup \halflineup \lineup \lineup;
					\draw \tlcoord{2.5}{1} \lineup;
					\draw \tlcoord{0}{1} \linedown;
					\draw \tlcoord{0}{1} \smalllineup \capright \smalllinedown \linedown;
					\draw \tlcoord{2.5}{1} \smalllinedown \cupright \smalllineup \lineup;
					\draw \tlcoord{0}{0} \maketlboxgreen{2}{$n, \varepsilon$};
					\draw \tlcoord{2.5}{0} \maketlboxgreen{2}{$n, \varepsilon$};
				\end{tikzpicture}
			} (n\geq2).
			\end{align*}
		\end{theo}
		
		Before the proof we establish diagrammatic properties of type $B$ projectors.
		
		\subsection{Diagrammatic properties}
		
		The following theorem gives a version of the property of Jones--Wenzl projectors: The annihilation of Temperley--Lieb generators.
		
		\begin{theo}[Type $B$ uniqueness] \label{Theorem on Existence of JW projectors of type B}
			Let $\eta\in \{{+},{-}\}$. The element $b_{n,\eta, \varepsilon}$ is the unique non-zero idempotent in $\TL(B_n)_{\varepsilon}$ which satisfies
			\begin{equation} \label{defining properties of positive type B JW} 
				s_0b_{n,\eta, \varepsilon}=\eta b_{n,\eta,\varepsilon}=b_{n,\eta,\varepsilon}s_0, \quad  U_ib_{n,\eta,\varepsilon}=0=b_{n,\eta,\varepsilon}U_i \, 
			\end{equation}
			for all $i=1,\ldots, n-1$.
		\end{theo}
		\begin{proof}
			The proof is analogous to the type $A$ proof, see for instance  \cite[Lem.~2]{lickorish1991}.
			For uniqueness one uses the decomposition
			\begin{equation} \label{nicuu direct sum decomposition}
				\TL(B_n)_{\varepsilon}={\shift{1, s_0}} \oplus (U_1,\ldots,U_{n-1}),
			\end{equation}
			The rest of the proof follows by simultaneous induction over the three statements:
			\begin{enumerate}[label = (\Roman*\textsubscript{$n$}), leftmargin=3\parindent]
				\item $b_{n \eta,\varepsilon}$ is an idempotent.
				\item $b_{n, \eta,\varepsilon}$ annihilates $U_j$, if $1\leq j \leq n-1$.
				\item $U_{n}b_{n, \eta,\varepsilon}U_n=-\varepsilon\frac{q^{n}+ q^{{-}n}}{q^{n-1}+ q^{{-}(n-1)}}U_{n}b_{n-1, \eta,\varepsilon}=-\varepsilon\frac{q^{n}+ q^{{-}n}}{q^{n-1}+ q^{{-}(n-1)}}b_{n-1, \eta,\varepsilon}U_{n}$. \qedhere
			\end{enumerate} 
		\end{proof}
	As consequence of \Cref{defi type B projectors}, \Cref{Theorem on Existence of JW projectors of type B} implies two stacking rules.
	\begin{koro}[Box absorption] \label{box absorption type B}
		Let $n\geq m\geq 1$, $\eta\in \{ {+}, {-}\}$. We have
		\begin{align*}
			b_{n,\eta,\varepsilon} b_{m,\eta,\varepsilon} &= b_{n,\eta,\varepsilon} =  b_{m,\eta, \varepsilon} b_{n,\eta, \varepsilon}
			\shortintertext{which we indicate diagrammatically by}
			\cbox{
				\begin{tikzpicture}[tldiagram, yscale=1/2, xscale=3/4]
					\draw \tlcoord{-2.5}{0} \lineup \lineup \lineup \lineup \lineup;
					\draw \tlcoord{-2.5}{1} \lineup \lineup \lineup \lineup \lineup;
					\draw \tlcoord{-2.5}{2} \lineup \lineup \lineup \lineup \lineup;
					\draw \tlcoord{1}{0} \maketlboxred{3}{$n, \varepsilon$};
					\draw \tlcoord{-1}{0} \maketlboxred{2}{$m, \varepsilon$};
					=
				\end{tikzpicture}
			}
			=
			\cbox{
				\begin{tikzpicture}[tldiagram, yscale=1/2, xscale=3/4]
					\draw \tlcoord{-2.5}{0} \lineup \lineup \lineup \lineup \lineup;
					\draw \tlcoord{-2.5}{1} \lineup \lineup \lineup \lineup \lineup;
					\draw \tlcoord{-2.5}{2} \lineup \lineup \lineup \lineup \lineup;
					\draw \tlcoord{0}{0} \maketlboxred{3}{$n, \varepsilon$};
				\end{tikzpicture}
			}
			\, &, \quad
			\cbox{
				\begin{tikzpicture}[tldiagram, yscale=1/2, xscale=3/4]
					\draw \tlcoord{-2.5}{0} \lineup \lineup \lineup \lineup \lineup;
					\draw \tlcoord{-2.5}{1} \lineup \lineup \lineup \lineup \lineup;
					\draw \tlcoord{-2.5}{2} \lineup \lineup \lineup \lineup \lineup;
					\draw \tlcoord{1}{0} \maketlboxblue{3}{$n,  \varepsilon$};
					\draw \tlcoord{-1}{0} \maketlboxblue{2}{$m,  \varepsilon$};
					=
				\end{tikzpicture}
			}
			=
			\cbox{
				\begin{tikzpicture}[tldiagram, yscale=1/2, xscale=3/4]
					\draw \tlcoord{-2.5}{0} \lineup \lineup \lineup \lineup \lineup;
					\draw \tlcoord{-2.5}{1} \lineup \lineup \lineup \lineup \lineup;
					\draw \tlcoord{-2.5}{2} \lineup \lineup \lineup \lineup \lineup;
					\draw \tlcoord{0}{0} \maketlboxblue{3}{$n, {-}$};
				\end{tikzpicture}
			} .
		\end{align*}
	\end{koro}
	
	\begin{koro}[Red-blue orthogonality] \label{red-blue orthogonality} 
		Let $n,m\geq 1$. The dominant and antidominant type $B$ projectors are orthogonal, i.e.\ we have
		\[
		\bp b_{m,{-},\varepsilon}=0=b_{m,{-},\varepsilon} \bp 
		\]
		which we indicate diagrammatically as
		\[
		\cbox{
			\begin{tikzpicture}[tldiagram, yscale=1/2, xscale=2/3]
				\draw \tlcoord{-2.5}{0} \lineup \lineup \lineup \lineup \lineup;
				\draw \tlcoord{-2.5}{1} \lineup \lineup \lineup \lineup \lineup;
				\draw \tlcoord{-2.5}{2} \lineup \lineup \lineup \lineup \lineup;
				\draw \tlcoord{1}{0} \maketlboxred{3}{$n, \varepsilon$};
				\draw \tlcoord{-1}{0} \maketlboxblue{2}{$m, \varepsilon$};
				=
			\end{tikzpicture}
		}
		\, = 0 = \, \cbox{
			\begin{tikzpicture}[tldiagram, yscale=1/2, xscale=2/3]
				\draw \tlcoord{-2.5}{0} \lineup \lineup \lineup \lineup \lineup;
				\draw \tlcoord{-2.5}{1} \lineup \lineup \lineup \lineup \lineup;
				\draw \tlcoord{-2.5}{2} \lineup \lineup \lineup \lineup \lineup;
				\draw \tlcoord{-1}{0} \maketlboxred{3}{$n, \varepsilon$};
				\draw \tlcoord{1}{0} \maketlboxblue{2}{$m, \varepsilon$};
				=
			\end{tikzpicture}
		} .
		\]
	\end{koro}
	These rules and \Cref{Theorem on Existence of JW projectors of type B} imply:
	\begin{koro}\label{type D properties}
		The element $d_{n,\varepsilon}$ is non-zero, idempotent and satisfies
		\[
		s_0d_{n,\varepsilon}s_0=d_{n,\varepsilon}, \quad d_{n,\varepsilon}U_i= 0 = U_i d_{n,\varepsilon} \, \text{ for all } i=0,1,\ldots, n-1.
		\]
	\end{koro}
	We conclude the following symmetry relating type $B$ and $D$ projectors.
	\begin{koro}[Symmetry] \label{symmetry}
		Let $\eta\in \{{+},{-}\}$. Then $b_{n,\eta,\varepsilon}=\frac{1}{2}(d_{n,\varepsilon}+\eta s_0d_{n,\varepsilon})$.
	\end{koro}
	\begin{proof}
		\Cref{type D properties} implies that the element $y_{\eta}\coloneqq\frac{1}{2}(d_{n,\varepsilon}+\eta s_0d_{n,\varepsilon})$ is idempotent, kills $U_1,\ldots, U_{n-1}$ and satisfies $s_0y_{\eta}=\eta y_{\eta}=y_{\eta}s_0$. Moreover it is non-zero, since $0\neq d_{n,\varepsilon}\in \TL(D_n)_{\varepsilon}$, which implies $d_{n,\varepsilon}$ and $s_0d_{n,\varepsilon}$ are linearly independent. Hence $y_{\eta}=b_{n,\eta,\varepsilon}$ by the uniqueness in \Cref{Theorem on Existence of JW projectors of type B}.
	\end{proof}
	\begin{koro}[Type $D_n$ uniqueness] \label{characterizing properties of type D projector}
			The projector $d_{n,\varepsilon}$ is the unique element in $\TL(D_n)_{\varepsilon}$, which satisfies the properties from \Cref{type D properties}.
	\end{koro}
	\begin{proof}
		Let $x\in \TL(D_n)_{\varepsilon}$ and assume it satisfies the properties from \Cref{type D properties}. By the argument as for \Cref{symmetry} we get $b_{n,\eta,\varepsilon}=\frac{1}{2}(x+\eta s_0x)$ for $\eta\in \{{+},{-}\}$. Hence by definition, 
			$d_{n,\varepsilon}=b_{n,{+},\varepsilon}+b_{n,{-},\varepsilon}=\frac{1}{2}(x+s_0x)+\frac{1}{2}(x- s_0x)=x. $\qedhere
	\end{proof}

		We can now also prove the recursive formula for type $D$ projectors.
		
		\begin{proof}[Proof of Theorem \ref{type D projector recursive description}]
			Assume $n\geq 2$. The recursive formula follows from
			\begin{align*}
				d_{n+1, \varepsilon}&=b_{n+1,{+}, \varepsilon}+b_{n+1,{-}, \varepsilon} \\
				&=\left(\bp + \lambda_n \bp U_n \bp\right) + \left(\bn +x_n \bn U_n \bn\right) \\
				&= d_{n, \varepsilon} +\lambda_n (\bp U_n \bp + \bn U_n \bn) \\
				&= d_{n, \varepsilon} + \lambda_n (\bp+\bn) U_n (\bp +\bn)  \\
				&= d_{n, \varepsilon} + \lambda_n d_{n+1, \varepsilon} U_n d_{n+1, \varepsilon} ,
			\end{align*}
			where $\lambda_n\coloneqq\varepsilon\frac{q^{n-1}+q^{-(n-1)}}{q^n+q^{-n}}$. We used in the non-trivial step that
			\[
			\bn U_n\bp = 0 = \bp U_n \bn.
			\]
			Indeed, \Cref{box absorption type B} implies that $b_{n,\eta, \varepsilon}=b_{1,\eta, \varepsilon}b_{n,\eta, \varepsilon}b_{1,\eta, \varepsilon}$ 
which we can apply to both $+,-$ to factor out a copy of $b_{1,{+}, \varepsilon}$ and $b_{1,{-}, \varepsilon}$ each. Since $b_{1,\eta, \varepsilon}$ and $U_n$ commute as $n\geq2$, the statement follows from \Cref{red-blue orthogonality} for $b_{1,+,\varepsilon}$ and $b_{1,-,\varepsilon}$.
		\end{proof}
		
		The recursive formulas for the type $D$ projectors in \Cref{defi type B projectors} appeared in non-quantized form in \cite{wedrich2018} as idempotents projecting onto the maximal and lowest weight space of the standard representation of $\sl2$. The next theorem lifts this interpretation to the quantum coideal setting. 
		
		Consider the tensor product decomposition of $V_{\varepsilon}^{\otimes n}$ over $\Vq$ in \Cref{Theorem decomposition of tensor power of V}. We write ${+}\coloneqq 1$ and ${-}\coloneqq -1$ to make notation compatible, so that we can consider the elements $b_{n,\varepsilon,\varepsilon}\coloneqq b_{n,\varepsilon}$ and $ b_{n,-\varepsilon,\varepsilon}\coloneqq b_{n,-\varepsilon} \in \TL(B_n)_{\varepsilon}$.
		
		\begin{theo} \label{theorem JW describe}
			Under the isomorphism $\Psi_{\varepsilon}\colon \TL(B_n)_{\varepsilon}\xrightarrow{\cong} \End_{\Vq}(V_\varepsilon^{\otimes n})$ from \Cref{Iso of TL and endo algebra}, the projector $b_{n,\varepsilon,\varepsilon}$ becomes the idempotent projection
			\[
			\pi_{n,\varepsilon,\varepsilon}\colon V_{\varepsilon}^{\otimes n}\rightarrow L(\varepsilon^{n-1}[n])\hookrightarrow V_{\varepsilon}^{\otimes n},
			\] and the projector $b_{n,-\varepsilon,\varepsilon}$ becomes the idempotent projection
			\[
			\pi_{n,-{\varepsilon},\varepsilon}\colon V_{\varepsilon}^{\otimes n}\rightarrow L(-\varepsilon^{n-1}[n])\hookrightarrow V_{\varepsilon}^{\otimes n}.
			\] 
			The type $D_n$ Jones--Wenzl idempotent $d_{n,\varepsilon}$ projects onto $L([n])\oplus L(-[n])$.
		\end{theo}
		\begin{proof}
			The strategy of the proof is to verify that the images $\Psi_{\varepsilon}^{-1}(\pi_{n,\varepsilon,\varepsilon})$, $\Psi_{\varepsilon}^{-1}(\pi_{n,-\varepsilon,\varepsilon})$ satisfy the characterizing properties of the projectors $b_{n,\varepsilon, \varepsilon}$ and $b_{n,-\varepsilon,\varepsilon}$ using the explicit vectors in the decomposition in \Cref{Theorem decomposition of tensor power of V}.
			First observe that $s_0v_{0,\varepsilon}=\varepsilon v_{0,\varepsilon}$ and similarly $s_0w_{0,\varepsilon}=-\varepsilon w_{0,\varepsilon}$.
			The subrepresentation $L(\varepsilon^{n-1}[n])$ is spanned by 
			\[
			v_{0,\varepsilon}\otimes v_{1,\varepsilon} \otimes \cdots \otimes v_{n-1,\varepsilon} \text{ for } \varepsilon=1 \quad v_{0,\varepsilon} \otimes w_{1, \varepsilon} \otimes v_{-2,\varepsilon} \otimes w_{3,\varepsilon} \otimes \cdots \text{ for } \varepsilon=-1
			\] 
			by \Cref{Theorem decomposition of tensor power of V}ii).
			In particular, $s_0$ acts on $L(\varepsilon^{n-1}[n])$ by $\varepsilon$ by the first observation. Similarly, the subrepresentation $L(-\varepsilon^{n-1}[n])$ is spanned by 
			\[
			w_{0,\varepsilon}\otimes w_{-1,\varepsilon}\otimes \cdots\otimes w_{-(n-1),\varepsilon} \text{ if } \varepsilon=1, \, w_{0,\varepsilon} \otimes v_{-1, \varepsilon} \otimes w_{2, \varepsilon} \otimes v_{-3, \varepsilon} \otimes \cdots \text{ if } \varepsilon=-1.
			\]
			Hence $s_0$ acts by $-\varepsilon$ on $L(-\varepsilon^{n-1}[n])$. We conclude that $\Psi^{-1}(\pi_{n,\varepsilon,\varepsilon})$ satisfies
			\[
			s_0\Psi^{-1}(\pi_{n,\varepsilon,\varepsilon})=\varepsilon\Psi^{-1}(\pi_{n,\varepsilon,\varepsilon})=\Psi^{-1}(\pi_{n,\varepsilon,\varepsilon})s_0,  
			\]
			and similarly that $\Psi^{-1}(\pi_{n,-{\varepsilon},\varepsilon})$ satisfies
			\[
			s_0\Psi^{-1}(\pi_{n,-\varepsilon,\varepsilon})=-\varepsilon\Psi^{-1}(\pi_{n,-\varepsilon,\varepsilon})=\Psi^{-1}(\pi_{n,\varepsilon,\varepsilon})s_0.
			\]
			Next in order to conclude that $U_i\Psi^{-1}(\pi_{n,\varepsilon,\varepsilon})=0$ and  $U_i\Psi^{-1}(\pi_{n,-\varepsilon,\varepsilon})=0$,
			observe $[n]$ and $-[n]$ do not appear as eigenvalues of $B$ acting on $V_{\varepsilon}^{\otimes l}$ by \Cref{Theorem decomposition of tensor power of V}~iv) whenever $l<n$. This shows that all morphisms $V_{\varepsilon}^{\otimes n}\rightarrow V_{\varepsilon}^{\otimes n}$ which factor through such $V_{\varepsilon}^{\otimes l}$ become zero, when composed with the projection onto $L([n])$ respectively $L(-[n])$. In particular this applies to the morphisms 
			\[
			\Psi_i(U_i)=\tlline \cdots \tlline \tlcupcap \tlline \cdots \tlline\colon V_{\varepsilon}^{\otimes n}\rightarrow V_{\varepsilon}^{\otimes (n-2)}\rightarrow V_{\varepsilon}^{\otimes n}.
			\] 
			This shows that $\Psi^{-1}(\pi_{n,\varepsilon,\varepsilon})$ (and  $\Psi^{-1}(\pi_{n,-\varepsilon,\varepsilon})$) satisfies the characterizing properties of $b_{n,\varepsilon,\varepsilon}$ (respectively $b_{n,-\varepsilon,\varepsilon}$), and hence they agree with these elements by  \Cref{Theorem on Existence of JW projectors of type B}. The claim for $d_{n,\varepsilon}$ follows from $d_{n,\varepsilon}=\bp+\bn=b_{n,\varepsilon,\varepsilon}+b_{n,-\varepsilon,\varepsilon}.$
		\end{proof}
		We finally illustrate how the diagrammatics can be used to compute $\Hom$-spaces.
		\begin{beis} \label{Counting dimensions via diagrams} \label{Counting dimensions via diagrams better}
			Let $a,b,c\geq 0$. Assume everything is type $1$ for simplicity. There are four kinds of $\Hom$-spaces over the coideal we are interested in:
			\begin{gather*}
			\Hom(V_{a}\otimes V_{b}, L([c])\oplus L([-c])), \quad \Hom(L([a])\oplus L([-a])\otimes V_{b}, V_{c}), \\
			\Hom(V_{a}\otimes V_{b}, V_{c}) , \quad 
			\Hom(L([a])\oplus L([-a])\otimes V_{b}, L([c])\oplus L([-c])).  
			\end{gather*}
			In order to determine these $\Hom$-spaces one can either use \Cref{general characters} or count all diagrams connecting three boxes of appropriate colors and numbers of strands. For example for $a=b=c=2$ one can see by \Cref{general characters}
			\[
				(L([2])\oplus L([-2]))\otimes V_2\cong L([4])\oplus L([2]) \oplus L([0])^2 \oplus L(-[2]) \oplus L(-[4]),
			\]
			which implies $\dim \Hom((L([2])\oplus L([-2]))\otimes V_2, L([2])\oplus L([-2]))=2$. The corresponding diagram basis consists of the two only possible connections
			\[
				\cbox
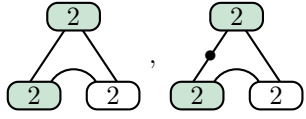
{
					\begin{tikzpicture}[tldiagram, scale=1/2]
						\draw \tlcoord{0}{0} -- ++ (4,6);
						\draw \tlcoord{0}{1} \xcapright{2};
						\draw \tlcoord{0}{4} -- ++ (-4,6);
						\draw \tlcoord{0}{0.25} \maketlboxgreen{1.5}{$2$};
						\draw \tlcoord{0}{3.25} \maketlboxnormal{1.5}{$2$};
						\draw \tlcoord{2}{1.75} \maketlboxgreen{1.5}{$2$};
					\end{tikzpicture}
				}, \cbox{
					\begin{tikzpicture}[tldiagram, scale=1/2]
						\draw \tlcoord{0}{0} -- ++ (2,3) \onedot -- ++ (2,3);
						\draw \tlcoord{0}{1} \xcapright{2};
						\draw \tlcoord{0}{4} -- ++ (-4,6);
						\draw \tlcoord{0}{0.25} \maketlboxgreen{1.5}{$2$};
						\draw \tlcoord{0}{3.25} \maketlboxnormal{1.5}{$2$};
						\draw \tlcoord{2}{1.75} \maketlboxgreen{1.5}{$2$};
					\end{tikzpicture}
				}
			\]
			between $d_2$, $a_1$ on the bottom and $d_2$ on the top.
			If one replaces the target $L([2])\oplus L([-2])$ with $V_2 \cong L([2])\oplus L(0)\oplus L([-2])$ over the coideal, the dimension of the $\Hom$-space increases by $2$. This corresponds to the two additional diagrams
			\begin{gather*}
				 \cbox
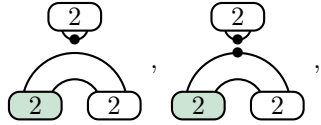
{
		\begin{tikzpicture}[tldiagram, scale=1/2]
			\draw \tlcoord{0}{0} \xcapright{4};
			\draw \tlcoord{0}{1} \xcapright{2};
			\draw \tlcoord{2}{1.5} \dcupright;
			\draw \tlcoord{0}{0.25} \maketlboxgreen{1.5}{$2$};
			\draw \tlcoord{0}{3.25} \maketlboxnormal{1.5}{$2$};
			\draw \tlcoord{2.25}{1.75} \maketlboxnormal{1.5}{$2$};
		\end{tikzpicture}
	}, \cbox{
	\begin{tikzpicture}[tldiagram, scale=1/2]
		\draw \tlcoord{0}{0} \dxcapright{4};
		\draw \tlcoord{0}{1} \xcapright{2};
		\draw \tlcoord{2}{1.5} \dcupright;
		\draw \tlcoord{0}{0.25} \maketlboxgreen{1.5}{$2$};
		\draw \tlcoord{0}{3.25} \maketlboxnormal{1.5}{$2$};
		\draw \tlcoord{2.25}{1.75} \maketlboxnormal{1.5}{$2$};
	\end{tikzpicture}
},
			\end{gather*}
			which become possible since the type $A$ projector doesn't kill a dotted cup. This is compatible with \Cref{beis decomposition}. Thinking further if one replaces the first tensor factor $L([2])\otimes L(-[2])$ in the source with $V_2=L([2])\oplus L(0)\oplus  L(-[2])$ one obtains
			\[
				V_2 \otimes V_2\cong L([4])\oplus L([2])^2 \oplus L([0])^3 \oplus L(-[2])^2 \oplus L(-[4])
			\]
			which increases the dimension of $\Hom$ by $3$. The additional three diagrams are:
			\begin{gather*}
			 \cbox
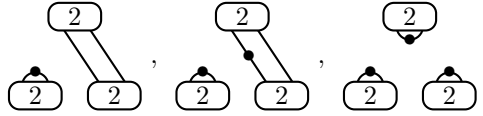
{
				\begin{tikzpicture}[tldiagram, scale=1/2]
					\draw \tlcoord{0.25}{0} \dcapright;
					\draw \tlcoord{0}{3.25} -- ++ (-2,3) -- ++ (-2,3);
					\draw \tlcoord{0}{4.25} -- ++ (-4,6);
					\draw \tlcoord{0}{0.25} \maketlboxnormal{1.5}{$2$};
					\draw \tlcoord{0}{3.25} \maketlboxnormal{1.5}{$2$};
					\draw \tlcoord{2}{1.75} \maketlboxnormal{1.5}{$2$};
				\end{tikzpicture}
			}, 	\, \cbox{
			\begin{tikzpicture}[tldiagram, scale=1/2]
				\draw \tlcoord{0.25}{0} \dcapright;
				\draw \tlcoord{0}{3.25} -- ++ (-2,3) \onedot -- ++ (-2,3);
				\draw \tlcoord{0}{4.25} -- ++ (-4,6);
				\draw \tlcoord{0}{0.25} \maketlboxnormal{1.5}{$2$};
				\draw \tlcoord{0}{3.25} \maketlboxnormal{1.5}{$2$};
				\draw \tlcoord{2}{1.75} \maketlboxnormal{1.5}{$2$};
			\end{tikzpicture}
		}, \, \cbox{
					\begin{tikzpicture}[tldiagram, scale=1/2]
						\draw \tlcoord{0.25}{0} \dcapright;
						\draw \tlcoord{0.25}{3} \dcapright;
						\draw \tlcoord{1.75}{1.5} \dcupright;
						\draw \tlcoord{0}{0.25} \maketlboxnormal{1.5}{$2$};
						\draw \tlcoord{0}{3.25} \maketlboxnormal{1.5}{$2$};
						\draw \tlcoord{2}{1.75} \maketlboxnormal{1.5}{$2$};
					\end{tikzpicture}
				}.
			\end{gather*}
		\end{beis}
		\section{Application: 
			\texorpdfstring{$\Theta$}{Theta}--networks} \label{section application to theta networks}
		We finish the paper by an outlook how our tools allow to define and compute $\Theta$-networks generalizing the networks arising in the Turaev--Viro invariants, \cite{flath95} and in the context of $3j$-symbols, see also \cite{stroppel2012}.
		First recall:
		\begin{defi}
			A triple $(a,b,c)\in \mN^{3}$ is called admissible if one (all) of the following equivalent conditions is satisfied:
			\begin{enumerate}
				\item $a\in \{b+c,b+c-2,b+c-4,\ldots, |b-c|\}$
				\item $b\in \{a+c,a+c-2,a+c-4,\ldots, |a-c|\}$
				\item $c\in \{a+b,a+b-2,a+b-4,\ldots, |a-b|\}$
			\end{enumerate}
			To any admissible triple $(a,b,c)$ we associate a triple $(i,j,k)\in  \mN^{3}$ by
			\[
			i\coloneqq \frac{b+c-a}{2}, \, j\coloneqq \frac{a+c-b}{2}, \, k\coloneqq \frac{a+b-c}{2}.
			\]
			In particular we can recover $(a,b,c)$ from $(i,j,k)$ as 
			$a=j+k, \quad b=i+k, \quad c=i+j$.
	
		\end{defi}
	\begin{lemm}
		The following are equivalent for a triple $(a,b,c)\in \mN^{3}$:\\
		\begin{minipage}[t]{6cm}
	\begin{enumerate}[labelindent=-6.8pt]
			\item[\textbullet]$(a,b,c)$ is admissible, \label{i1}
			\item[\textbullet] $\Hom_{\Uq}(V_a\otimes V_b, V_c)\neq 0$, \label{i2}
			\item[\textbullet] $\Hom_{\Vq}(V_a\otimes V_b, V_c)\neq 0$,  \label{i3}
			\item[\textbullet] $\Hom_{\Vq}(V_a\otimes V_b, L([c]))\neq 0$, \label{i6} 
			\item[\textbullet] $\Hom_{\Vq}(V_a\otimes V_b, L([-c]))\neq 0$, \label{i7}
			\end{enumerate}	
		\end{minipage}
		\begin{minipage}[t]{6.4cm}	
			\begin{enumerate}[labelindent=-6pt]
			\item[\textbullet] $\Hom_{\Vq}(L([a])\otimes V_b, V_c)\neq 0$,  \label{i4} 
			\item[\textbullet] $\Hom_{\Vq}(L([-a])\otimes V_b, V_c)\neq 0$,  \label{i5}
			\item[\textbullet] $\Hom_{\Vq}(L([a])\otimes V_b, L([c]))\neq 0$, \label{i8} 
			\item[\textbullet] $\Hom_{\Vq}(L([-a])\otimes V_b, L([-c]))\neq 0$. \label{i9} 
			\end{enumerate}
			\end{minipage}
	\end{lemm}
	\begin{proof}
		This follows from \Cref{general characters}.
	\end{proof}
\begin{rema}
	By the Clebsch--Gordan rule the $\Hom$-space in \ref{i2} is at most one dimensional. The up-to-scalar unique non-zero morphism is a unique type $A$ planar diagram connecting the type $A$ projectors $a_{a-1}, a_{b-1}, a_{c-1}$, where $i$ many strands connect $a_{b-1}$ and $a_{c-1}$, $j$ many connect $a_{a-1}$ and $a_{c-1}$, and $k$ many connect $a_{a-1}$ and $a_{b-1}$. 
	As we saw in \Cref{Counting dimensions via diagrams better} the $\Hom$-spaces \ref{i3}- \ref{i9} are in general higher dimensional with bases labeled by all possible type $B$ diagrams one can draw between them. Hence in order to define analogues $\Theta$-networks one has to not only fix an admissible triple, but also fix a type $B$ diagram connecting the projectors. For the rest of the section we work with type $D$ projectors, instead of type $B$, and fix the standard type $A$ connection. The first choice does not lose any information, since the type $D$ projector carries all of the information about the type $B$ projectors. A full treatment of all the choices is material of another paper, which also treats $3j$-symbols similarly to \cite[Theorem 8]{stroppel2012} with respect to the bases we defined.
\end{rema}
	
		\begin{defi}
			Let $\varepsilon\in\{1,-1\}$. Let $(a,b,c)$ be an admissible triple. A  $\Theta$-network
			$\Theta_{a,b,c,\varepsilon}^D$ of type $D$ is the rational function $f(1)$, where $f\colon\C(q)\rightarrow \C(q)$ is given by
		\[
			\Theta_{a,b,c,\varepsilon}^D \coloneqq \cbox
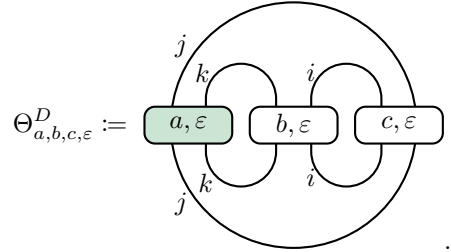
{
				\begin{tikzpicture}[tldiagram, scale=2/3]
					\draw \tlcoord{0}{0}    \xcapright{7} \xcupleft{7};
					\draw \tlcoord{-0.5}{1}   \lineup \xcapright{2} \linedown \xcupleft{2};
					\draw \tlcoord{-0.5}{4}   \lineup \xcapright{2} \linedown \xcupleft{2};
					\draw \tlcoord{0}{0} \maketlboxgreen{2}{$a, \varepsilon$};
					\draw \tlcoord{0}{3} \maketlboxnormal{2}{$b, \varepsilon$};
					\draw \tlcoord{0}{6} \maketlboxnormal{2}{$c, \varepsilon$};
					\node (n1) at \tlcoord{1.5}{0.3} {$j$};
					\node (n1) at \tlcoord{1}{0.9} {$k$};
					\node (n1) at \tlcoord{1}{4} {$i$};
					\node (n1) at \tlcoord{-1.5}{0.3} {$j$};
					\node (n1) at \tlcoord{-1.1}{1} {$k$};
					\node (n1) at \tlcoord{-1}{4} {$i$};
				\end{tikzpicture}.	
			}
			\]
		The outer cup and cap should be read as $j$ nested strands, the left small cup/cap as $k$ nested strands and the right small cup/cap as $i$ nested strands.
		\end{defi}
		The following theorem gives a recursive description of type $D$ $3j$-symbols.
		\begin{theo} \label{recursive formula for type D Theta networks}
			For all $a\geq 1$ we have
			\begin{enumerate}
				\item $\Theta_{a,0,a,\varepsilon}^D=(-\varepsilon)^a \left(q^a+q^{-a}\right)$
				\item $\Theta_{a,1,a-1,\varepsilon}^D=(-\varepsilon)^a \left(q^a+q^{-a}\right)$ 
				\item $\Theta_{a-1,1,a,\varepsilon}^D=(-\varepsilon)^a \frac{[a]}{[a-1]}\left(q^{a-1}+q^{-(a-1)}\right)$
			\end{enumerate}
			Moreover for any admissible triple $(a,b,c)$ with $a\geq 1, b\geq 2, c\geq 1$ the type $D$ $3j$-symbol $\Theta_{a,b,c,\varepsilon}^D$ satisfies the following recursion: 
			\[
			\Theta_{a,b,c,\varepsilon}^D=-\varepsilon\frac{[c]}{[c-1]}\Theta_{a,b-1,c-1,\varepsilon}^D+\varepsilon\frac{[(a+b-c)/2]^2}{[b][b-1]}\Theta_{a,b-2,c,\varepsilon} 
			\]	
		\end{theo}
		Before we outline a proof of the theorem, we recall some useful formulas for computing $\Theta$-networks in type $A$ \cite[Lemma 1,2]{masbaum94}. Both lemmas describe rules how to delete strands, which is helpful for recursion.
		\begin{lemm} \label{lemma to move strand to right}
			Let $i,j,k\in \mN^3$ such that $i,j\geq 1$. Then we have:
			\[
			\cbox{
				\begin{tikzpicture}[tldiagram, scale=9/10]
					\draw \tlcoord{0}{0} \linewave{2}{2};
					\draw \tlcoord{0}{1}  \xcapright{3.5};
					\draw \tlcoord{-0.3}{1.7} \xcupright{1} \smalllineup \xcapright{1};
					\draw \tlcoord{0}{5.5} \linewave{2}{-2};
					\draw \tlcoord{0}{0} \maketlboxnormal{2.5}{$i+j-1,\varepsilon$};
					\draw \tlcoord{0}{4} \maketlboxnormal{2.5}{$j+k,\varepsilon$};
					\draw \tlcoord{2}{2} \maketlboxnormal{2.5}{$i+k,\varepsilon$};
					\node (n1) at \tlcoord{1.2}{1} {$i$};
					\node (n1) at \tlcoord{1.3}{2.5} {$j-1$};
					\node (n1) at \tlcoord{1.2}{4.6} {$k$};
					\node (n1) at \tlcoord{-0.5}{3} {$1$};
				\end{tikzpicture}	
			} \! \! \! \! \!=  \varepsilon^{j-1}\frac{[i]}{[i+j-1]}\cbox{
				\begin{tikzpicture}[tldiagram, scale=9/10]
					\draw \tlcoord{0}{0} \linewave{2}{2};
					\draw \tlcoord{0}{1} \xcapright{3.5};
					\draw \tlcoord{0}{5.5} \linewave{2}{-2};
					\draw \tlcoord{0}{0} \maketlboxnormal{2.5}{$i+j-2,\varepsilon$};
					\draw \tlcoord{0}{4} \maketlboxnormal{2.5}{$j+k,\varepsilon$};
					\draw \tlcoord{2}{2} \maketlboxnormal{2.5}{$i+k,\varepsilon$};
					\node (n1) at \tlcoord{1.2}{0.5} {$i-1$};
					\node (n1) at \tlcoord{1.3}{2.5} {$j-1$};
					\node (n1) at \tlcoord{1.2}{5} {$k+1$};
				\end{tikzpicture}	
			}
			\]
		\end{lemm}
		\begin{lemm} \label{Lemma closing up one strand in type A}
			Let $n\geq 1$. Then we have
			\[	
			\cbox{
				\begin{tikzpicture}[tldiagram, yscale=1/2, xscale=2/3, scale=-1]
					\draw \tlcoord{-1.5}{0} \lineup \lineup \lineup ;
					\draw \tlcoord{-1.5}{1} \lineup \lineup \lineup ;
					\draw \tlcoord{-0.5}{2} \lineup   \xcapright{2}   \linedown  \xcupleft{2};
					\draw \tlcoord{0}{0} \maketlboxnormal{3}{$n, \varepsilon$};
				\end{tikzpicture}
			} = -\varepsilon \frac{[n]}{[n-1]} \cbox{
				\begin{tikzpicture}[tldiagram, yscale=1/2]
					\draw \tlcoord{-1.5}{0} \lineup \lineup \lineup ;
					\draw \tlcoord{-1.5}{1} \lineup \lineup \lineup ;
					\draw \tlcoord{0}{0} \maketlboxnormal{2}{$n-1, \varepsilon$};
				\end{tikzpicture}
			}.
			\]
		\end{lemm}
		
		The characterizing property of the type $D$ Jones--Wenzl projector directly implies: 
		
		\begin{lemm} \label{lemma color absorption}
			We have $d_{n,\varepsilon}a_{n-1,\varepsilon}=d_{n,\varepsilon}=a_{n-1,\varepsilon}d_{n,\varepsilon}$. Diagrammatically speaking:
			\[
			\cbox{
				\begin{tikzpicture}[tldiagram, yscale=1/2, xscale=2/3]
					\draw \tlcoord{-2.5}{0} \lineup \lineup \lineup \lineup \lineup;
					\draw \tlcoord{-2.5}{1} \lineup \lineup \lineup \lineup \lineup;
					\draw \tlcoord{-2.5}{2} \lineup \lineup \lineup \lineup \lineup;
					\draw \tlcoord{1}{0} \maketlboxgreen{3}{$n, \varepsilon$};
					\draw \tlcoord{-1}{0} \maketlboxnormal{3}{$n, \varepsilon$};
					=
				\end{tikzpicture}
			}
			\,=\,
			\cbox{
				\begin{tikzpicture}[tldiagram, yscale=1/2, xscale=2/3]
					\draw \tlcoord{-2.5}{0} \lineup \lineup \lineup \lineup \lineup;
					\draw \tlcoord{-2.5}{1} \lineup \lineup \lineup \lineup \lineup;
					\draw \tlcoord{-2.5}{2} \lineup \lineup \lineup \lineup \lineup;
					\draw \tlcoord{0}{0} \maketlboxgreen{3}{$n, \varepsilon$};
				\end{tikzpicture}
			} 
			\, = \, \cbox{
				\begin{tikzpicture}[tldiagram, yscale=1/2, xscale=2/3]
					\draw \tlcoord{-2.5}{0} \lineup \lineup \lineup \lineup \lineup;
					\draw \tlcoord{-2.5}{1} \lineup \lineup \lineup \lineup \lineup;
					\draw \tlcoord{-2.5}{2} \lineup \lineup \lineup \lineup \lineup;
					\draw \tlcoord{1}{0} \maketlboxnormal{3}{$n, \varepsilon$};
					\draw \tlcoord{-1}{0} \maketlboxgreen{3}{$n, \varepsilon$};
					=
				\end{tikzpicture}
			},
			\]
		\end{lemm}

		The next lemma is a type $D$ version of Lemma \ref{Lemma closing up one strand in type A}.
		\begin{lemm} \label{lemma closing up one strand}
			Let $n\geq 2$. We have the following formula:
			\[
			U_nd_{n, \varepsilon}U_n=-\varepsilon \frac{q^{n}+q^{-n}}{q^{n-1}+q^{-(n-1)}} U_nd_{n-1,\varepsilon}
			\]
			Diagramatically this becomes:
			\[
			\cbox{
				\begin{tikzpicture}[tldiagram, yscale=1/2, xscale=2/3]
					\draw \tlcoord{-1.5}{0} \lineup \lineup \lineup ;
					\draw \tlcoord{-1.5}{1} \lineup \lineup \lineup ;
					\draw \tlcoord{-0.5}{2} \lineup   \xcapright{2}   \linedown  \xcupleft{2};
					\draw \tlcoord{0}{0} \maketlboxgreen{3}{$n, \varepsilon$};
				\end{tikzpicture}
			} = -\varepsilon \frac{q^{n}+q^{-n}}{q^{n-1}+q^{-(n-1)}} \cbox{
				\begin{tikzpicture}[tldiagram, yscale=1/2]
					\draw \tlcoord{-1.5}{0} \lineup \lineup \lineup ;
					\draw \tlcoord{-1.5}{1} \lineup \lineup \lineup ;
					\draw \tlcoord{0}{0} \maketlboxgreen{2}{$n-1, \varepsilon$};
				\end{tikzpicture}
			}
			\]
		\end{lemm}
		\begin{proof}
			We already proved this statement for $\bp$ and $\bn$ as part III$_n$ of the proof of Theorem \ref{Theorem on Existence of JW projectors of type B}. Now use $d_{n,\varepsilon}=\bp+\bn$. The diagrammatic version follows from surjectivity of $\cap\colon V_{\varepsilon}\otimes V_{\varepsilon}\to \C(q)$ and injectivity of $\cup\colon \C(q) \to V_{\varepsilon}\otimes V_{\varepsilon}$. 
		\end{proof}
		As direct consequence of the previous lemma we obtain:
		\begin{koro} \label{Corollary dimension formula}
			Let $n\geq 1$. We have the following trace formula: 
			\[
			\cbox{
				\begin{tikzpicture}[tldiagram, yscale=5/7, scale=4/7]
					\draw \tlcoord{-0.5}{0}   \lineup \xcapright{6}   \linedown \xcupleft{6};
					\draw \tlcoord{-0.5}{1}   \lineup \xcapright{4}   \linedown \xcupleft{4};
					\draw \tlcoord{-0.5}{2}   \lineup \xcapright{2} \linedown \xcupleft{2};
					\draw \tlcoord{0}{0} \maketlboxgreen{3}{$n, \varepsilon$};
				\end{tikzpicture}
			} = (-\varepsilon)^{n} \left({q^{n}+q^{-n}}\right)
			\]
		\end{koro}
		\begin{proof}
			This statement is clear for $n=1,2$ by the explicit descriptions. For the induction step use Lemma \ref{lemma closing up one strand}.  
		\end{proof}
		Now we can prove of the main theorem:
		\begin{proof}[Proof of Theorem \ref{recursive formula for type D Theta networks}]
			The formulas for $\Theta_{a,0,a,\varepsilon}^D$, $\Theta_{a,0,a-1,\varepsilon}^D$ and $\Theta_{a-1,0,a,\varepsilon}^D$ are direct consequences of Lemma \ref{Lemma closing up one strand in type A}, Lemma \ref{lemma color absorption}, Lemma \ref{lemma closing up one strand} and Corollary \ref{Corollary dimension formula}.
			For the recursion formula we sketch the (diagrammatic) proof:
			First apply the type $A$ Wenzl recursion formula (Definition \ref{recursive formula type A}) to the middle Jones--Wenzl projector on $b$ strands. This yields two summands. For the first one use Lemma \ref{Lemma closing up one strand in type A} to obtain $-\varepsilon\frac{[c]}{[c-1]}\Theta_{a,b-1,c-1,\varepsilon}^D$.
			For the second summand use Lemma \ref{lemma color absorption} twice to write the type $D$ Jones--Wenzl projector $d_{n,\varepsilon}$ as $a_{n-1,\varepsilon}d_{n,\varepsilon}a_{n-1,\varepsilon}$.
			Then bend the two $A_{n-1}$ projectors over to the right. Afterwards apply Lemma \ref{lemma to move strand to right} to the right hand side twice. Then use that $a_{n-1,\varepsilon}$ is idempotent twice and bend the type $A$ Jones--Wenzl projectors to the left to obtain $\varepsilon\frac{[(a+b-c)/2]^2}{[b][b-1]}\Theta_{a,b-2,c,\varepsilon}$ as the second summand.
		\end{proof}
		\bibliographystyle{alpha}
\bibliography{literature}
\end{document}